\documentclass[final,onefignum,onetabnum]{siamart190516}

\usepackage{lipsum}
\usepackage{amsfonts}
\usepackage{graphicx}
\usepackage{epstopdf}
\usepackage{algorithmic}
\usepackage{arydshln}
\usepackage{enumitem}

\ifpdf
  \DeclareGraphicsExtensions{.eps,.pdf,.png,.jpg}
\else
  \DeclareGraphicsExtensions{.eps}
\fi

\newcommand{\define}{\stackrel{\text{\tiny def}}{=}}

\newsiamremark{remark}{Remark}
\newsiamremark{hypothesis}{Hypothesis}
\crefname{hypothesis}{Hypothesis}{Hypotheses}
\newsiamthm{claim}{Claim}

\headers{General Tail Bounds for Random Tensors Summation}{Shih Yu Chang}

\title{General Tail Bounds for Random Tensors Summation: Majorization Approach}

\author{Shih Yu Chang\thanks{Department of Applied Data Science, San Jose State University, San Jose, CA 95192, USA 
  (\email{ shihyu.chang@sjsu.edu } ).} \and
Yimin Wei\thanks{School of Mathematical Sciences, Fudan University, Shanghai, China (\email{ymwei@fudan.edu.cn }).}
}

\usepackage{amsopn}

\ifpdf
\hypersetup{
  pdftitle={General Tail Bounds for Random Tensors Summation: Majorization Approach},
  pdfauthor={Shih Yu Chang}
}
\fi

\begin{document}

\maketitle

\begin{abstract}
In recent years, tensors have been applied to different applications in science and engineering fields. In order to establish theory about tail bounds of the tensors summation behavior, this work extends previous work by considering the tensors summation tail behavior of the top $k$-largest singular values of a function of the tensors summation, instead of the largest/smallest singular value of the tensors summation directly (identity function) explored in~\cite{chang2020convenient}. Majorization and antisymmetric tensor product tools are main techniques utilized to establish inequalities for unitarily norms of multivariate tensors. The Laplace transform method is integrated with these inequalities for unitarily norms of multivariate tensors to give us tail bounds estimation for Ky Fan $k$-norm for a function of the tensors summation. By restricting different random tensor conditions, we obtain generalized tensor Chernoff and Bernstein inequalities. 
\end{abstract}



\section{Introduction}

\subsection{From Random Matrices to Radom Tensors}\label{sec:From Random Matrices to Radom Tensors}



A random matrix is a matrix-valued random variable—that is, a matrix in which some or all entries are random variables. Random matrices have played an important role in numerical linear algebra~\cite{MR41539}, quantum mecaniques~\cite{MR730191}, neural networks~\cite{wainrib2013topological}, communication theory~\cite{tulino2004random}, robust control~\cite{hosoe2017robust}, etc. Many important properties of scientific and engineering systems can be modelled by matrix formulations. In order to consider a high-dimensional system, it is often more convenient to consider tensors, or \emph{multidimensional data}, instead of matrices (two-dimensional data). 

In recent years, tensors have been applied to different applications in science and engineering~\cite{qi2017tensor}. In data processing fields, tensor theory applications include unsupervised separation of unknown mixtures of data signals~\cite{wu2010robust, mirsamadi2016generalized}, signals filtering~\cite{muti2007survey}, network signal processing~\cite{shen2020topology, shen2017tensor, fu2015joint} and image processing~\cite{ko2020fast, jiang2020framelet}. In wireless communication applications, tensors are applied to model high-dimensional communication channels, e.g., MIMO (multi-input multi-output) code-division~\cite{de2008constrained, zhijin2018blind}, radar communications~\cite{nion2010tensor, sidiropoulos2000parallel}. In numerical multilinear algebra computations, tensors can be applied to solve multilinear system of equations~\cite{wang2019neural}, high-dimensional data fitting/regression~\cite{MR3395816}, tensor complementary problem~\cite{MR3947912}, tensor eigenvalue problem~\cite{MR3479021}, etc. In machine learning, tensors are also used to characterize data with \emph{coupling effects}, for example, tensor decomposition methods have been reported recently to establish the latent-variable models, such as topic models in~\cite{anandkumar2015tensor}, and the method of moments for undertaking the \emph{Latent Dirichlet Allocation} (LDA) in~\cite{sidiropoulos2017tensor}. Nevertheless, all these applications assume that systems modelled by tensors are fixed and such assumption is not true and practical in problems involving tensor formulations. In recent years, there are more works beginning to develop theory about random tensors, see~\cite{MR3616422,MR3783911,MR4140540}, and references therein. In this work, we will apply majorization techniques to establish inequalities for unitarily norms of multivariate tensors and these inequalities will be used to derive more general tensor Chernoff and Bernstein inequalities.

\subsection{Technical Results}\label{sec:Technical Results} 


Majorization is an effective tool for proving norm and trace inequalities of linear operators, see~\cite{MR2759813} for more advanced topics therein. An important topic investigated by majorization technique is to answer the following question:
\\
\textit{\textbf{Question 1}: Given a continuous function $f$, two matrices $\mathbf{C}$ and $\mathbf{D}$, what is the relationship between the majorization order among singular values of $f(\mathbf{C})$ and $f(\mathbf{D})$, and the matrix norm relation of $f(\mathbf{C})$ and $f(\mathbf{D})$?}
\\

In Prop. 4.4.13 in~\cite{hiai2010matrix}, they showed that the singular values of $\mathbf{C}$ are weakly majorized by the singular values of $\mathbf{D}$ if and only if $\left\Vert \mathbf{C} \right\Vert_{\rho} \leq \left\Vert \mathbf{D} \right\Vert_{\rho}$ for every unitarily invariant norm~\footnote{The exact unitarily invariant norm definition will be provided in Section~\ref{sec:Unitarily Invariant Tensor Norms}.} $\left\Vert \cdot \right\Vert_{\rho}$, e.g., Schatten $p$-norms and trace norm. In this previous work, they only find this relationship in \textbf{Question 1} when the function $f$ is an identity map. In another recent work~\cite{hiai2017generalized}, they apply log-majorization to find relationship in \textbf{Question 1} by considering more general continous function $f$ and utilize this new relationship to prove multivariate generalizations of the Araki–Lieb–Thirring inequality and the Golden–Thompson inequality for unitarily invariant norm of matrices.

In this work, we generalize this approach to answer \textbf{Question 1} in tensors settings. Our first main result is to prove multivaraite tensor norm inequalities given by Theorem~\ref{thm:Multivaraite Tensor Norm Inequalities intro}. Besides extending from matrices to tensors, our work tries to simplify some steps in the approach adopated by~\cite{hiai2017generalized}, for example,  Lemma~\ref{lma:15} is simplified from Lemma 15 in~\cite{hiai2017generalized}. 


\begin{theorem}\label{thm:Multivaraite Tensor Norm Inequalities intro}
Let $\mathcal{C}_i \in \mathbb{C}^{I_1 \times \cdots \times I_N \times I_1 \times \cdots \times I_N}$ be positive definite Hermitian tensors for $1 \leq i \leq n$ with Hermitian rank $r$, $\left\Vert \cdot \right\Vert_{\rho}$ be a unitarily invaraint norm with corresponding gauge function $\rho$. For any continous function $f:(0, \infty) \rightarrow [0, \infty)$ such that $x \rightarrow \log f(e^x)$ is convex on $\mathbb{R}$, we have 
\begin{eqnarray}\label{eq1:thm:Multivaraite Tensor Norm Inequalities intro}
\left\Vert  f \left( \exp \left( \sum\limits_{i=1}^n \log \mathcal{C}_i\right)   \right)  \right\Vert_{\rho} &\leq& \exp \int_{- \infty}^{\infty} \log \left\Vert f \left( \left\vert \prod\limits_{i=1}^{n}  \mathcal{C}_i^{1 + \iota t} \right\vert\right)\right\Vert_{\rho} \beta_0(t) dt ,
\end{eqnarray}
where $\iota= \sqrt{-1}$ and $\beta_0(t) = \frac{\pi}{2 (\cosh (\pi t) + 1)}$.

For any continous function $g(0, \infty) \rightarrow [0, \infty)$ such that $x \rightarrow g (e^x)$ is convex on $\mathbb{R}$, we have 
\begin{eqnarray}\label{eq2:thm:Multivaraite Tensor Norm Inequalities intro}
\left\Vert  g \left( \exp \left( \sum\limits_{i=1}^n \log \mathcal{C}_i\right)   \right)  \right\Vert_{\rho} &\leq& \int_{- \infty}^{\infty} \left\Vert g \left( \left\vert \prod\limits_{i=1}^{n}  \mathcal{C}_i^{1 + \iota t} \right\vert\right)\right\Vert_{\rho} \beta_0(t) dt.
\end{eqnarray}
\end{theorem}

There are two main technical tools required by this work to build those tensor probability bounds. The first is \emph{Laplace transform method}, which provides a systematic way to give tail bounds for the sum of scalar random variables. In~\cite{MR1966716}, the authors apply Laplace transform method to bound the largest eigenvalue with the matrix setting, i.e., the tail probability for the maximum eigenvalue of the sum of Hermitian matrices is controlled by a matrix version of the moment-generating function. They prove following:
\begin{eqnarray}\label{eq:Ahlswede Winter bound for max eigen}
\mathrm{Pr} \left( \lambda_{\max}\left( \sum\limits_{i} \mathbf{X}_i \right) \geq \theta    \right) \leq \inf\limits_{t > 0} \Big\{ e^{-t \theta} \mathbb{E} \mathrm{Tr} \exp \left( \sum\limits_{i} t \mathbf{X}_i \right)\Big\}.
\end{eqnarray}
In this work, we extend the Laplace transform method to tensors and utilize Theorem~\ref{thm:Multivaraite Tensor Norm Inequalities intro} to obtain the following Ky Fan $k$-norm bounds for the tail behavior of a function of tensors summation. 


\begin{theorem}\label{thm:Ky Fan norm prob bound for fun of tensors sum intro}
Consider a sequence $\{ \mathcal{X}_j  \in \mathbb{C}^{I_1 \times \cdots \times I_N  \times I_1 \times \cdots \times I_N} \}$ of independent, random, Hermitian tensors. Let $g$ be a polynomial function with degree $n$ and nonnegative coeffecients $a_0, a_1, \cdots, a_n$ raised by power $s \geq 1$, i.e., $g(x) = \left(a_0 + a_1 x  +\cdots + a_n x^n \right)^s$. Suppose following condition is satisfied:
\begin{eqnarray}\label{eq:special cond intro}
g \left( \exp\left(t \sum\limits_{j=1}^{m} \mathcal{X}_j \right)\right)  \succeq \exp\left(t g \left( \sum\limits_{j=1}^{m} \mathcal{X}_j   \right) \right)~~\mbox{almost surely},
\end{eqnarray}
where $t > 0$. Then, we have 
\begin{eqnarray}\label{eq1:thm:Ky Fan norm prob bound for fun of tensors sum intro}
\mathrm{Pr} \left( \left\Vert g\left( \sum\limits_{j=1}^{m} \mathcal{X}_j  \right)\right\Vert_{(k)}  \geq \theta \right) \leq ~~~~~~~~~~~~~~~~~~~~~~~~~~~~~~~~~~~~~~~~~~~~ \nonumber \\
(n+1)^{s-1}\inf\limits_{t, p_j} e^{- \theta t }\left(k a_0^s + \sum\limits_{l=1}^{n}a^{ls}_l \sum\limits_{j=1}^m \frac{ \mathbb{E} \left\Vert \exp\left( p_j  l s t \mathcal{X}_j \right) \right\Vert_{(k)} }{p_j}     \right).
\end{eqnarray}
where $\sum\limits_{j=1}^m \frac{1}{p_j} =1$ and $p_j > 0$. 
\end{theorem}

If we restrict following conditions: $\mathcal{X}_i \succeq \mathcal{O}$ and $\lambda_{\max}(\mathcal{X}_i) \leq \mathrm{R}, \mbox{~~ almost surely;}$ to random tensors $\mathcal{X}_i$, we can further bound the term of $\mathbb{E} \left\Vert \exp\left( p_j  l s t \mathcal{X}_j \right) \right\Vert_{(k)}$ to get the following generalized tensor Chernoff bound provided by Theorem~\ref{thm:Generalized Tensor Chernoff Bound intro}. 


\begin{theorem}[Generalized Tensor Chernoff Bound]\label{thm:Generalized Tensor Chernoff Bound intro}
Consider a sequence $\{ \mathcal{X}_j  \in \mathbb{C}^{I_1 \times \cdots \times I_N  \times I_1 \times \cdots \times I_N} \}$ of independent, random, Hermitian tensors. Let $g$ be a polynomial function with degree $n$ and nonnegative coeffecients $a_0, a_1, \cdots, a_n$ raised by power $s \geq 1$, i.e., $g(x) = \left(a_0 + a_1 x  +\cdots + a_n x^n \right)^s$ with $s \geq 1$. Suppose following condition is satisfied:
\begin{eqnarray}\label{eq:special cond Chernoff Bound intro}
g \left( \exp\left(t \sum\limits_{j=1}^{m} \mathcal{X}_j \right)\right)  \succeq \exp\left(t g \left( \sum\limits_{j=1}^{m} \mathcal{X}_j   \right) \right)~~\mbox{almost surely},
\end{eqnarray}
where $t > 0$. Moreover, we require  
\begin{eqnarray}
\mathcal{X}_i \succeq \mathcal{O} \mbox{~~and~~} \lambda_{\max}(\mathcal{X}_i) \leq \mathrm{R}
\mbox{~~ almost surely.}
\end{eqnarray}
Then we have following inequality:
\begin{eqnarray}\label{eq1:thm:GeneralizedTensorChernoffBound intro}
\mathrm{Pr} \left( \left\Vert g\left( \sum\limits_{j=1}^{m} \mathcal{X}_j  \right)\right\Vert_{(k)}  \geq \theta \right)  \leq  (n+1)^{s-1} \inf\limits_{t > 0} e^{- \theta t } \cdot ~~~~~~~~~~~~~~~~~~~~~~~~~~~~~~ \nonumber \\
 \left\{ ka_0^s + \sum\limits_{l=1}^{n} \sum\limits_{j=1}^m \frac{k a_l^{ls}}{m} \left[ 1 +\left( e^{mlsRt} - 1 \right) \overline{\sigma_1(\mathcal{X}_j)}  +  C \left( e^{mlsRt} - 1 \right) \Xi(\mathcal{X}_j) \right] \right\},
\end{eqnarray}
where $C$ is a constant and $\Xi(\mathcal{X}_j)$ is determined from the expectation of entries from the tensor $\mathcal{X}_j$ defined by Eq.~\eqref{eq:abbr of sigma 1 of a random tensor}.
\end{theorem}

On the other hand, if we consider following conditions to random tensors $\mathcal{X}_i$:
\begin{eqnarray}
\mathbb{E} \mathcal{X}_j = \mathcal{O} \mbox{~~and~~} \mathcal{X}^p_j \preceq \frac{p! \mathcal{A}_j^2}{2}
\mbox{~~ almost surely for $p=2,3,4,\cdots$;}
\end{eqnarray}
we will get the following generalized tensor Bernstein bound provided by Theorem~\ref{thm:Generalized Tensor Bernstein Bound intro}.


\begin{theorem}[Generalized Tensor Bernstein Bound]\label{thm:Generalized Tensor Bernstein Bound intro}
Consider a sequence $\{ \mathcal{X}_j  \in \mathbb{C}^{I_1 \times \cdots \times I_N  \times I_1 \times \cdots \times I_N} \}$ of independent, random, Hermitian tensors. Let $g$ be a polynomial function with degree $n$ and nonnegative coeffecients $a_0, a_1, \cdots, a_n$ raised by power $s \geq 1$, i.e., $g(x) = \left(a_0 + a_1 x  +\cdots + a_n x^n \right)^s$ with $s \geq 1$. Suppose following condition is satisfied:
\begin{eqnarray}\label{eq:special cond Bernstein Bound intro}
g \left( \exp\left(t \sum\limits_{j=1}^{m} \mathcal{X}_j \right)\right)  \succeq \exp\left(t g \left( \sum\limits_{j=1}^{m} \mathcal{X}_j   \right) \right)~~\mbox{almost surely},
\end{eqnarray}
where $t > 0$, and we also have 
\begin{eqnarray}
\mathbb{E} \mathcal{X}_j = \mathcal{O} \mbox{~~and~~} \mathcal{X}^p_j \preceq \frac{p! \mathcal{A}_j^2}{2}
\mbox{~~ almost surely for $p=2,3,4,\cdots$.}
\end{eqnarray}
Then we have following inequality:
\begin{eqnarray}\label{eq1:thm:Generalized Tensor Bernstein Bound intro}
\mathrm{Pr} \left( \left\Vert g\left( \sum\limits_{j=1}^{m} \mathcal{X}_j  \right)\right\Vert_{(k)}  \geq \theta \right) & \leq & (n+1)^{s-1} \inf\limits_{t > 0} e^{- \theta t }k \cdot \nonumber \\
&  & \left\{a_0^s + \sum\limits_{l=1}^{n} \sum\limits_{j=1}^m  a_l^{ls} \left[ \frac{1}{m} + \frac{ m (lst)^2 \sigma_1(\mathcal{A}_j^2) }{2 (1 - m lst)}+ lstC \Upsilon(\mathcal{X}_j) \right] \right\},
\end{eqnarray}
where $C$ is a constant and $\Upsilon(\mathcal{X}_j)$ is determined from the expectation of entries from the tensor $\mathcal{X}_j$ defined by Eq.~\eqref{eq:abbr of sigma 1 of a random tensor zero mean}.
\end{theorem}

The bounds provided by Theorems~\ref{thm:Ky Fan norm prob bound for fun of tensors sum intro},~\ref{thm:Generalized Tensor Chernoff Bound intro} and~\ref{thm:Generalized Tensor Bernstein Bound intro} are more general than those random tensor bounds given by~\cite{chang2020convenient} since we can consider the concentration behavior of the top $k$-largest singular values, instead of the largest singular value, and a function of tensors summation, instead of the identity map of tensors summation. We have another work based on similar majorization techniques to consider Chernoff expander bounds for random tensors without independent random tensors assumptions~\cite{chang2021tensor}.

\subsection{Paper Organization}\label{sec:Paper Organization}

The paper is organized as follows. Preliminaries of tensors and basic majorization notations are given in Section~\ref{sec:Fundamentals of Tensors and Majorization}. In Section~\ref{sec:Multivariate Tensor Norm Inequalities}, we will develop theorems about majorization and log majorization with integral average which will be applied to derive bounds for unitarily norms of multivariate tensors. Compared to work~\cite{chang2020convenient}, more generalized tensor Chernoff and Bernstein inequalities are discussed in Section~\ref{sec:New Tensor Inequalities}. Finally, the conclusions are given in Section~\cref{sec:Conclusions}.

\section{Fundamentals of Tensors and Majorization}\label{sec:Fundamentals of Tensors and Majorization}


\subsection{Tensors Preliminaries}\label{sec:Tensors Preliminaries}

Throughout this work, scalars are represented by lower-case letters (e.g., $d$, $e$, $f$, $\ldots$), vectors by boldfaced lower-case letters (e.g., $\mathbf{d}$, $\mathbf{e}$, $\mathbf{f}$, $\ldots$), matrices by boldfaced capitalized letters (e.g., $\mathbf{D}$, $\mathbf{E}$, $\mathbf{F}$, $\ldots$), and tensors by calligraphic letters (e.g., $\mathcal{D}$, $\mathcal{E}$, $\mathcal{F}$, $\ldots$), respectively. Tensors are multiarrays of values which are higher-dimensional generalizations from vectors and matrices. Given a positive integer $N$, let $[N] \define \{1, 2, \cdots ,N\}$. An \emph{order-$N$ tensor} (or \emph{$N$-th order tensor}) denoted by $\mathcal{X} \define (x_{i_1, i_2, \cdots, i_N})$, where $1 \leq i_j = 1, 2, \ldots, I_j$ for $j \in [N]$, is a multidimensional array containing $\prod_{n=1}^N I_n$ entries. 
Let $\mathbb{C}^{I_1 \times \cdots \times I_N}$ and $\mathbb{R}^{I_1 \times \cdots \times I_N}$ be the sets of the order-$N$ $I_1 \times \cdots \times I_N$ tensors over the complex field $\mathbb{C}$ and the real field $\mathbb{R}$, respectively. For example, $\mathcal{X} \in \mathbb{C}^{I_1 \times \cdots \times I_N}$ is an order-$N$ multiarray, where the first, second, ..., and $N$-th dimensions have $I_1$, $I_2$, $\ldots$, and $I_N$ entries, respectively. Thus, each entry of $\mathcal{X}$ can be represented by $x_{i_1, \cdots, i_N}$. For example, when $N = 3$, $\mathcal{X} \in \mathbb{C}^{I_1 \times I_2 \times I_3}$ is a third-order tensor containing entries $x_{i_1, i_2, i_3}$'s.

Without loss of generality, one can partition the dimensions of a tensor into two groups, say $M$ and $N$ dimensions, separately. Thus, for two order-($M$+$N$) tensors: $\mathcal{X} \define (x_{i_1, \cdots, i_M, j_1, \cdots,j_N}) \in \mathbb{C}^{I_1 \times \cdots \times I_M\times
J_1 \times \cdots \times J_N}$ and $\mathcal{Y} \define (y_{i_1, \cdots, i_M, j_1, \cdots,j_N}) \in \mathbb{C}^{I_1 \times \cdots \times I_M\times
J_1 \times \cdots \times J_N}$, according to~\cite{MR3913666}, the \emph{tensor addition} $\mathcal{X} + \mathcal{Y}\in \mathbb{C}^{I_1 \times \cdots \times I_M\times
J_1 \times \cdots \times J_N}$ is given by 
\begin{eqnarray}\label{eq: tensor addition definition}
(\mathcal{X} + \mathcal{Y} )_{i_1, \cdots, i_M, j_1 \times \cdots \times j_N} &\define&
x_{i_1, \cdots, i_M, j_1 \times \cdots \times j_N} \nonumber \\
& &+ y_{i_1, \cdots, i_M, j_1 \times \cdots \times j_N}. 
\end{eqnarray}
On the other hand, for tensors $\mathcal{X} \define (x_{i_1, \cdots, i_M, j_1, \cdots,j_N}) \in \mathbb{C}^{I_1 \times \cdots \times I_M\times
J_1 \times \cdots \times J_N}$ and $\mathcal{Y} \define (y_{j_1, \cdots, j_N, k_1, \cdots,k_L}) \in \mathbb{C}^{J_1 \times \cdots \times J_N\times K_1 \times \cdots \times K_L}$, according to~\cite{MR3913666}, the \emph{Einstein product} (or simply referred to as \emph{tensor product} in this work) $\mathcal{X} \star_{N} \mathcal{Y} \in  \mathbb{C}^{I_1 \times \cdots \times I_M\times
K_1 \times \cdots \times K_L}$ is given by 
\begin{eqnarray}\label{eq: Einstein product definition}
\lefteqn{(\mathcal{X} \star_{N} \mathcal{Y} )_{i_1, \cdots, i_M,k_1 \times \cdots \times k_L} \define} \nonumber \\ &&\sum\limits_{j_1, \cdots, j_N} x_{i_1, \cdots, i_M, j_1, \cdots,j_N}y_{j_1, \cdots, j_N, k_1, \cdots,k_L}. 
\end{eqnarray}
Note that we will often abbreviate a tensor product $\mathcal{X} \star_{N} \mathcal{Y}$ to ``$\mathcal{X} \hspace{0.05cm}\mathcal{Y}$'' for notational simplicity in the rest of the paper. 
This tensor product will be reduced to the standard matrix multiplication as $L$ $=$ $M$ $=$ $N$ $=$ $1$. Other simplified situations can also be extended as tensor–vector product ($M >1$, $N=1$, and $L=0$) and tensor–matrix product ($M>1$ and $N=L=1$). In analogy to matrix analysis, we define some basic tensors and elementary tensor operations as follows. 

\begin{definition}\label{def: zero tensor}
A tensor whose entries are all zero is called a \emph{zero tensor}, denoted by $\mathcal{O}$. 
\end{definition}

\begin{definition}\label{def: identity tensor}
An \emph{identity tensor} $\mathcal{I} \in  \mathbb{C}^{I_1 \times \cdots \times I_N\times
J_1 \times \cdots \times J_N}$ is defined by 
\begin{eqnarray}\label{eq: identity tensor definition}
(\mathcal{I})_{i_1 \times \cdots \times i_N\times
j_1 \times \cdots \times j_N} \define \prod_{k = 1}^{N} \delta_{i_k, j_k},
\end{eqnarray}
where $\delta_{i_k, j_k} \define 1$ if $i_k  = j_k$; otherwise $\delta_{i_k, j_k} \define 0$.
\end{definition}

In order to define \emph{Hermitian} tensor, the \emph{conjugate transpose operation} (or \emph{Hermitian adjoint}) of a tensor is specified as follows.  
\begin{definition}\label{def: tensor conjugate transpose}
Given a tensor $\mathcal{X} \define (x_{i_1, \cdots, i_M, j_1, \cdots,j_N}) \in \mathbb{C}^{I_1 \times \cdots \times I_M\times J_1 \times \cdots \times J_N}$, its conjugate transpose, denoted by
$\mathcal{X}^{H}$, is defined by
\begin{eqnarray}\label{eq:tensor conjugate transpose definition}
(\mathcal{X}^H)_{ j_1, \cdots,j_N,i_1, \cdots, i_M}  \define  
 x^*_{i_1, \cdots, i_M,j_1, \cdots,j_N},
\end{eqnarray}
where the star $*$ symbol indicates the complex conjugate of the complex number $x_{i_1, \cdots, i_M,j_1, \cdots,j_N}$. If a tensor $\mathcal{X}$ satisfies $ \mathcal{X}^H = \mathcal{X}$, then $\mathcal{X}$ is a \emph{Hermitian tensor}. 
\end{definition}
We will use symbol $\iota$ to represent $\sqrt{-1}$.

Following definition is about untiary tensors.
\begin{definition}\label{def: unitary tensor}
Given a tensor $\mathcal{U} \define (u_{i_1, \cdots, i_M, i_1, \cdots,i_M}) \in \mathbb{C}^{I_1 \times \cdots \times I_M\times I_1 \times \cdots \times I_M}$, if
\begin{eqnarray}\label{eq:unitary tensor definition}
\mathcal{U}^H \star_M \mathcal{U} = \mathcal{U} \star_M \mathcal{U}^H = \mathcal{I} \in \mathbb{C}^{I_1 \times \cdots \times I_M\times I_1 \times \cdots \times I_M},
\end{eqnarray}
then $\mathcal{U}$ is a \emph{unitary tensor}. 
\end{definition}
In this work, the symbol $\mathcal{U}$ is reserved for a unitary tensor. 

\begin{definition}\label{def: inverse of a tensor}
Given a \emph{square tensor} $\mathcal{X} \define (x_{i_1, \cdots, i_M, j_1, \cdots,j_M}) \in \mathbb{C}^{I_1 \times \cdots \times I_M\times I_1 \times \cdots \times I_M}$, if there exists $\mathcal{X} \in \mathbb{C}^{I_1 \times \cdots \times I_M\times I_1 \times \cdots \times I_M}$ such that 
\begin{eqnarray}\label{eq:tensor invertible definition}
\mathcal{X} \star_M \mathcal{X} = \mathcal{X} \star_M \mathcal{X} = \mathcal{I},
\end{eqnarray}
then $\mathcal{X}$ is the \emph{inverse} of $\mathcal{X}$. We usually write $\mathcal{X} \define \mathcal{X}^{-1}$ thereby. 
\end{definition}

We also list other crucial tensor operations here. The \emph{trace} of a square tensor is equivalent to the summation of all diagonal entries such that 
\begin{eqnarray}\label{eq: tensor trace def}
\mathrm{Tr}(\mathcal{X}) \define \sum\limits_{1 \leq i_j \leq I_j,\hspace{0.05cm}j \in [M]} \mathcal{X}_{i_1, \cdots, i_M,i_1, \cdots, i_M}.
\end{eqnarray}
The \emph{inner product} of two tensors $\mathcal{X}$, $\mathcal{Y} \in \mathbb{C}^{I_1 \times \cdots \times I_M\times J_1 \times \cdots \times J_N}$ is given by 
\begin{eqnarray}\label{eq: tensor inner product def}
\langle \mathcal{X}, \mathcal{Y} \rangle \define \mathrm{Tr}\left(\mathcal{X}^H \star_M \mathcal{Y}\right).
\end{eqnarray}
According to Eq.~\eqref{eq: tensor inner product def}, the \emph{Frobenius norm} of a tensor $\mathcal{X}$ is defined by 
\begin{eqnarray}\label{eq:Frobenius norm}
\left\Vert \mathcal{X} \right\Vert \define \sqrt{\langle \mathcal{X}, \mathcal{X} \rangle}.
\end{eqnarray}

%
%

From Theorem 5.2 in~\cite{ni2019hermitian}, every Hermitian tensor $\mathcal{H} \in  \mathbb{C}^{I_1 \times \cdots \times I_N \times I_1 \times \cdots \times I_N}$ has following decomposition
\begin{eqnarray}\label{eq:Hermitian Eigen Decom}
\mathcal{H} &=& \sum\limits_{i=1}^r \lambda_i \mathcal{U}_i \otimes \mathcal{U}^{H}_i, \mbox{
~with~~$\langle \mathcal{U}_i, \mathcal{U}_i \rangle =1$ and $\langle \mathcal{U}_i, \mathcal{U}_j \rangle = 0$ for $i \neq j$,}
\end{eqnarray}
where $\lambda_i \in \mathbb{R}$ and $\otimes$ denotes for Kronecker product. The values $\lambda_i$ are named as \emph{Hermitian eigevalues}, and the minimum integer of $r$ to decompose a Hermitian tensor as in Eq.~\eqref{eq:Hermitian Eigen Decom} is called \emph{Hermitian tensor rank}. A \emph{postive Hermitian tensor} is a Hermitian tensor with all \emph{Hermitian eigevalues} are positive. A \emph{nonnegative Hermitian tensor} is a Hermitian tensor with all \emph{Hermitian eigevalues} are nonegative.  The \emph{Hermitian determinant}, denoted as $\det\nolimits_H(\mathcal{A})$,  is defined as the product of $\lambda_i$ of the tensor $\mathcal{A}$. 

\textbf{Axioms For Hermitian Determinnant}
The Hermitian determinant is a mapping from Hermitian tensor to a real number satisfying following: 
\begin{enumerate}
	\item $\det_H (\mathcal{I}) = 1$.
	\item  $\det_H (A \star B ) = \det_H(A) \det_H(B)$.
\end{enumerate}

Then, we have $\left\vert \det(U) \right\vert = 1$, $\det(A) = \det(\Lambda)$ (the product of all eigenvalues)

\subsection{Unitarily Invariant Tensor Norms}\label{sec:Unitarily Invariant Tensor Norms}


Let us represent the Hermitian eigenvalues of a Hermitian tensor $\mathcal{H} \in \mathbb{C}^{I_1 \times \cdots \times I_N \times I_1 \times \cdots \times I_N} $ in decreasing order by the vector $\vec{\lambda}(\mathcal{H}) = (\lambda_1(\mathcal{H}), \cdots, \lambda_r(\mathcal{H}))$, where $r$ is the Hermitian rank of the tensor $\mathcal{H}$. We use $\mathbb{R}_{\geq 0} (\mathbb{R}_{> 0})$ to represent a set of nonnegative (positive) real numbers. Let $\left\Vert \cdot \right\Vert_{\rho}$ be a unitarily invariant tensor norm, i.e., $\left\Vert \mathcal{H}\star_N \mathcal{U}\right\Vert_{\rho} = \left\Vert \mathcal{U}\star_N \mathcal{H}\right\Vert_{\rho} = \left\Vert \mathcal{H}\right\Vert_{\rho} $,  where $\mathcal{U}$ is any unitary tensor. Let $\rho : \mathbb{R}_{\geq 0}^r \rightarrow \mathbb{R}_{\geq 0}$ be the corresponding gauge function that satisfies H$\ddot{o}$lder’s inequality so that 
\begin{eqnarray}\label{eq:def gauge func and general unitarily invariant norm}
\left\Vert \mathcal{H} \right\Vert_{\rho} = \left\Vert |\mathcal{H}| \right\Vert_{\rho} = \rho(\vec{\lambda}( | \mathcal{H} | ) ),
\end{eqnarray}
where $ |\mathcal{H}|  \define \sqrt{\mathcal{H}^H \star_N \mathcal{H}} $. The bijective correspondence between symmetric gauge functions on $\mathbb{R}_{\geq 0}^r$ and unitarily invariant norms is due to von Neumann~\cite{fan1955some}. 

Several popular norms can be treated as special cases of unitarily invariant tensor norm. The first one is Ky Fan like $k$-norm~\cite{fan1955some} for tensors. For $k \in \{1,2,\cdots, r \}$, the Ky Fan $k$-norm~\cite{fan1955some} for tensors  $\mathcal{H}  \mathbb{C}^{I_1 \times \cdots \times I_N \times I_1 \times \cdots \times I_N} $, denoted as $\left\Vert \mathcal{H}\right\Vert_{(k)}$, is defined as:
\begin{eqnarray}\label{eq: Ky Fan k norm for tensors}
\left\Vert \mathcal{H}\right\Vert_{(k)} \define \sum\limits_{i=1}^{k} \lambda_i(  |\mathcal{H}|  ).
\end{eqnarray}
If $k=1$,  the Ky Fan $k$-norm for tensors is the tensor operator norm, denoted as $ \left\Vert \mathcal{H} \right\Vert$. The second one is Schatten $p$-norm for tensors, denoted as $\left\Vert \mathcal{H}\right\Vert_{p}$, is defined as:
\begin{eqnarray}\label{eq: Schatten p norm for tensors}
\left\Vert \mathcal{H}\right\Vert_{p} \define (\mathrm{Tr}|\mathcal{H}|^p )^{\frac{1}{p}},
\end{eqnarray}
where $ p \geq 1$. If $p=1$, it is the trace norm. The third one is $k$-trace norm, denoted as $\mathrm{Tr}_k[\mathcal{H}]$, defined by ~\cite{huang2020generalizing}. It is 
\begin{eqnarray}\label{eq: de k-trace norm for tensors}
\mathrm{Tr}_k[\mathcal{H}] \define \sum\limits_{1 \leq i_1 < i_2 < \cdots i_k \leq r} \lambda_{i_1} \lambda_{i_1}  \cdots \lambda_{i_k} 
\end{eqnarray}
where $ 1 \leq k \leq r$. If $k=1$, $\mathrm{Tr}_k[\mathcal{H}]$ is reduced as trace norm. 

Following inequality is the extension of H\"{o}lder inequality to gauge function $\rho$ which will be used by later to prove majorization relations. 
\begin{lemma}\label{lma:Holder inquality for gauge function}
For $n$ nonnegative real vectors with the dimension $r$, i.e., $\mathbf{b}_i = (b_{i_1}, \cdots, b_{i_r}) \in \mathbb{R}_{\geq 0}^r$, and $\alpha > 0$ with $\sum\limits_{i=1}^n \alpha_i = 1$, we have 
\begin{eqnarray}\label{eq1:lma:Holder inquality for gauge function}
\rho\left( \prod\limits_{i=1}^n b_{i_1}^{\alpha_i},  \prod\limits_{i=1}^n b_{i_2}^{\alpha_i}, \cdots,  \prod\limits_{i=1}^n b_{i_r}^{\alpha_i}  \right) \leq  \prod\limits_{i=1}^n \rho(\mathbf{b}_i)^{\alpha_i} 
\end{eqnarray}
\end{lemma}
\textbf{Proof:}
This proof is based on mathematical induction. The base case for $n=2$ has been shown by Theorem IV.1.6 from~\cite{bhatia2013matrix}. 

We assume that Eq.~\eqref{eq1:lma:Holder inquality for gauge function} is true for $n=m$, where $m > 2$. Let $\odot$ be the component-wise product (Hadamard product) between two vectors.  Then, we have 
\begin{eqnarray}\label{eq2:lma:Holder inquality for gauge function}
\rho\left( \prod\limits_{i=1}^{m+1} b_{i_1}^{\alpha_i},  \prod\limits_{i=1}^{m+1} b_{i_2}^{\alpha_i}, \cdots,  \prod\limits_{i=1}^{m+1} b_{i_r}^{\alpha_i}  \right) = 
\rho\left( \odot_{i=1}^{m+1} \mathbf{b}_i^{\alpha_i}  \right),
\end{eqnarray}
where $\odot_{i=1}^{m+1} \mathbf{b}_i^{\alpha_i}$ is defined as $\left( \prod\limits_{i=1}^{m+1} b_{i_1}^{\alpha_i},  \prod\limits_{i=1}^{m+1} b_{i_2}^{\alpha_i}, \cdots,  \prod\limits_{i=1}^{m+1} b_{i_r}^{\alpha_i}  \right)$ with $\mathbf{b}_i^{\alpha_i} \define (b_{i_1}^{\alpha_i}, \cdots, b_{i_r}^{\alpha_i})$. Under such notations, Eq.~\eqref{eq2:lma:Holder inquality for gauge function} can be bounded as  
\begin{eqnarray}\label{eq3:lma:Holder inquality for gauge function}
\rho\left( \odot_{i=1}^{m+1} \mathbf{b}_i^{\alpha_i}  \right) &= &
\rho\left( \left( \odot_{i=1}^{m} \mathbf{b}_i^{\frac{\alpha_i}{ \sum\limits_{j=1}^m \alpha_j }  } \right)^{\sum\limits_{j=1}^m \alpha_j } \odot \mathbf{b}_{m+1}^{\alpha_{m+1}}\right) \nonumber \\
& \leq & \left[ \rho^{\sum\limits_{j=1}^m \alpha_j } \left( \odot_{i=1}^{m} \mathbf{b}_i^{\frac{\alpha_i}{ \sum\limits_{j=1}^m \alpha_j }  }  \right)  \right] \cdot \rho( \mathbf{b}_{m+1})^{\alpha_{m+1}} \leq  \prod\limits_{i=1}^{m+1} \rho(\mathbf{b}_i)^{\alpha_i}. 
\end{eqnarray}
By mathematical induction, this lemma is proved. $\hfill \Box$

\subsection{Antisymmetric Tensor Product}\label{sec:Antisymmetric Tensor Product}

Let $\mathfrak{H}$ be a Hilbert space of dimension $r$, $\mathfrak{L}(\mathfrak{H})$ be the set of tensors (linear operators) on $\mathfrak{H}$. Two tensors $\mathcal{A}, \mathcal{B} \in \mathfrak{L}(\mathfrak{H})$ is said $\mathcal{A} \geq \mathcal{B}$ if $\mathcal{A} - \mathcal{B} $ is a nonnegative Hermitian tensor. For any $k \in \{1,2,\cdots,r\}$, let $\mathfrak{H}^{\otimes k}$ be the $k$-th tensor power of the space $\mathfrak{H}$ and let $\mathfrak{H}^{\wedge k}$ be the antisymmetric subspace of $\mathfrak{H}^{\otimes k}$. We define function $\wedge^k : \mathfrak{L}( \mathfrak{H}) \rightarrow \mathfrak{L}(\mathfrak{H}^{\wedge k})$ as mapping any tensor $\mathcal{A}$ to the restriction of $\mathcal{A}^{\otimes k} \in \mathfrak{L}( \mathfrak{H}^{\otimes k} )$ to the antisymmetric subspace $\mathfrak{H}^{\wedge k}$ of $\mathfrak{H}^{\otimes k}$. Following lemma summarizes several useful properties of such antisymmetric tensor products.

\begin{lemma}\label{lma:antisymmetric tensor product properties}
Let $\mathcal{A}, \mathcal{B}, \mathcal{C} \in \mathbb{C}^{I_1 \times \cdots \times I_N \times I_1 \times \cdots \times I_N}$ be tensors in $\mathfrak{L}(\mathfrak{H})$, and $\mathcal{D} \in \mathbb{C}^{I_1 \times \cdots \times I_N \times I_1 \times \cdots \times I_N} $ be Hermitian tensors from $\mathfrak{H}$ with Hermitian rank $r$. For any $k \in \{1,2,\cdots,r\}$, we have 
\begin{enumerate}[label={[\arabic*]}]
	\item $(\mathcal{A}^{\wedge k})^H = (\mathcal{A}^H)^{\wedge k}$,.
	\item $(\mathcal{A}^{\wedge k}) \star_N (\mathcal{B}^{\wedge k})= (\mathcal{A}\star_N \mathcal{B})^{\wedge k}$. 
	\item If $\lim\limits_{i \rightarrow \infty} \left\Vert \mathcal{A}_i -  \mathcal{A} \right\Vert \rightarrow 0$ , then $\lim\limits_{i \rightarrow \infty} \left\Vert \mathcal{A}^{\wedge k}_i -  \mathcal{A}^{\wedge k} \right\Vert \rightarrow 0$.
	\item If $\mathcal{C} \geq \mathcal{O}$ (zero tensor), then $\mathcal{C}^{\wedge k} \geq \mathcal{O}$ and $(\mathcal{C}^p)^{\wedge k} = (\mathcal{C}^{\wedge k})^p$ for all $p \in \mathbb{R}_{> 0 }$.
    \item  $|\mathcal{A}|^{\wedge k} = | \mathcal{A}^{\wedge k}|$.
    \item  If $\mathcal{D} \geq \mathcal{O}$ and $\mathcal{D}$ is invertibale,  $(\mathcal{D}^z)^{\wedge k} = (\mathcal{D}^{\wedge k})^z$ for all $z \in \mathbb{D}$.
    \item  $\left\Vert \mathcal{A}^{\wedge k} \right\Vert = \prod\limits_{i=1}^{k} \lambda_i ( | \mathcal{A} |)$.
\end{enumerate}
\end{lemma}
\textbf{Proof:}
Facts $\textit{[1]}$ and $\textit{[2]}$ are the restrictions of the associated relations $(\mathcal{A}^H)^{\otimes k} = (\mathcal{A}^{\otimes k})^H$ and $(\mathcal{A} \star_N \mathcal{B})^{\otimes k} = (\mathcal{A}^{\otimes k})\star_N (\mathcal{B}^{\otimes k})$ to $\mathfrak{H}^{\wedge k}$. The fact $\textit{[3]}$ is true since, if $\lim\limits_{i \rightarrow \infty} \left\Vert \mathcal{A}_i -  \mathcal{A} \right\Vert \rightarrow 0$, we have $\lim\limits_{i \rightarrow \infty} \left\Vert \mathcal{A}^{\otimes k}_i -  \mathcal{A}^{\otimes k} \right\Vert \rightarrow 0$ and the asscoaited restrictions of $\mathcal{A}_i^{\otimes k}, \mathcal{A}^{\otimes k}$ to the antisymmetric subspace $\mathfrak{H}^k$. 

For the fact $\textit{[4]}$, if $\mathcal{C} \geq \mathcal{O}$, then we have $\mathcal{C}^{\wedge k} = ((\mathcal{C}^{1/2})^{\wedge k})^H \star_N  ((\mathcal{C}^{1/2})^{\wedge k}) \geq   \mathcal{O}$ from facts  $\textit{[1]}$ and $\textit{[2]}$. If $p$ is ratonal, we have  $(\mathcal{C}^p)^{\wedge k} = (\mathcal{C}^{\wedge k})^p$  from the fact $\textit{[2]}$, and the equality $(\mathcal{C}^p)^{\wedge k} = (\mathcal{C}^{\wedge k})^p$ is also true for any $p > 0$ if we apply the fact $\textit{[3]}$ to approximate any irrelational numbers by rational numbers.

Because we have 
\begin{eqnarray}
|\mathcal{A}|^{\wedge k} =  \left( \sqrt{\mathcal{A}^H \mathcal{A}}\right)^{\wedge k}  =   \sqrt{ (\mathcal{A}^{\wedge k})^H \mathcal{A}^{\wedge k}  }=  | \mathcal{A}^{\wedge k}|,
\end{eqnarray}
from facts $\textit{[1]}$, $\textit{[2]}$ and $\textit{[4]}$, so the fact $\textit{[5]}$ is valid. 

For the fact $\textit{[6]}$, if $z  < 0$, the fact $\textit{[6]}$ is true for all $z \in \mathbb{R}$ by applying the fact $\textit{[4]}$ to $\mathcal{D}^{-1}$. Since we can apply the definition $\mathcal{D}^{z} \define \exp(z \ln \mathcal{D})$ to have
\begin{eqnarray}
\mathcal{C}^p &=& \mathcal{D}^z~~\leftrightarrow~~\mathcal{C} = \exp\left(\frac{z}{p} \ln \mathcal{D} \right),
\end{eqnarray}
where $\mathcal{C} \geq \mathcal{O}$. The general case of any $z \in \mathbb{C}$ is also true by applying the fact $\textit{[4]}$ to $\mathcal{C} = \exp(\frac{z}{p} \ln \mathcal{D})$. 

For the fact $\textit{[7]}$ proof, it is enough to prove the case that $\mathcal{A} \geq \mathcal{O}$ due to the fact $\textit{[5]}$. Then, there exisits a set of orthogonal tensors $\{\mathcal{U}_1, \cdots, \mathcal{U}_r\}$ such that $\mathcal{A} \star_N \mathcal{U}_i = \lambda_i   \mathcal{U}_i$ for $1 \leq i \leq r$. We then have 
\begin{eqnarray}
\mathcal{A}^{\wedge k} \left(\mathcal{U}_{i_1} \wedge \cdots \wedge \mathcal{U}_{i_k}\right)
&=& \mathcal{A}\star_N  \mathcal{U}_{i_1} \wedge \cdots \wedge \mathcal{A}\star_N  \mathcal{U}_{i_k}  \nonumber \\
&=& \left( \prod\limits_{i=1}^{k} \lambda_i ( | \mathcal{A} |)  \right) \mathcal{U}_{i_1} \wedge \cdots \wedge \mathcal{U}_{i_k},
\end{eqnarray}
where $1 \leq i_1 < i_2 < \cdots < i_k \leq r$. Hence, $\left\Vert \mathcal{A}^{\wedge k} \right\Vert = \prod\limits_{i=1}^{k} \lambda_i ( | \mathcal{A} |)$. $\hfill \Box$

\subsection{Majorization}\label{sec:Majorization}


In this subsection, we will discuss majorization and several lemmas about majorization which will be used at later proofs. 

Let $\mathbf{x} = [x_1, \cdots,x_r] \in \mathbb{R}^{r}, \mathbf{y} = [y_1, \cdots,y_r] \in \mathbb{R}^{r}$ be two vectors with following orders among entries $x_1 \geq \cdots \geq x_r$ and $y_1 \geq \cdots \geq y_r$, \emph{weak majorization} between vectors $\mathbf{x}, \mathbf{y}$, represented by $\mathbf{x} \prec_{w} \mathbf{y}$, requires following relation for  vectors $\mathbf{x}, \mathbf{y}$:
\begin{eqnarray}\label{eq:weak majorization def}
\sum\limits_{i=1}^k x_i \leq \sum\limits_{i=1}^k y_i,
\end{eqnarray}
where $k \in \{1,2,\cdots,r\}$. \emph{Majorization} between vectors $\mathbf{x}, \mathbf{y}$, indicated by $\mathbf{x} \prec \mathbf{y}$, requires following relation for vectors $\mathbf{x}, \mathbf{y}$:
\begin{eqnarray}\label{eq:majorization def}
\sum\limits_{i=1}^k x_i &\leq& \sum\limits_{i=1}^k y_i,~~\mbox{for $1 \leq k < r$;} \nonumber \\
\sum\limits_{i=1}^r x_i &=& \sum\limits_{i=1}^r y_i,~~\mbox{for $k = r$.}
\end{eqnarray}

For $\mathbf{x}, \mathbf{y} \in \mathbb{R}^r_{\geq 0}$ such that  $x_1 \geq \cdots \geq x_r$ and $y_1 \geq \cdots \geq y_r$,  \emph{weak log majorization} between vectors $\mathbf{x}, \mathbf{y}$, represented by $\mathbf{x} \prec_{w \log} \mathbf{y}$, requires following relation for vectors $\mathbf{x}, \mathbf{y}$:
\begin{eqnarray}\label{eq:weak log majorization def}
\prod\limits_{i=1}^k x_i \leq \prod\limits_{i=1}^k y_i,
\end{eqnarray}
where $k \in \{1,2,\cdots,r\}$, and \emph{log majorization} between vectors $\mathbf{x}, \mathbf{y}$, represented by $\mathbf{x} \prec_{\log} \mathbf{y}$, requires equality for $k=r$ in Eq.~\eqref{eq:weak log majorization def}. If $f$ is a single variable function, $f(\mathbf{x})$ represents a vector of $[f(x_1),\cdots,f(x_r)]$. From Lemma 1 in~\cite{hiai2017generalized}, we have 
\begin{lemma}\label{lma:Lemma 1 Gen Log Hiai}
(1) For any convex function $f: [0, \infty) \rightarrow [0, \infty)$, if we have $\mathbf{x} \prec \mathbf{y}$, then $f(\mathbf{x}) \prec_{w} f(\mathbf{y})$. \\
(2) For any convex function and non-decreasing $f: [0, \infty) \rightarrow [0, \infty)$, if we have $\mathbf{x} \prec_{w} \mathbf{y}$, then $f(\mathbf{x}) \prec_{w} f(\mathbf{y})$. \\
\end{lemma}

Another lemma is from Lemma 12 in~\cite{hiai2017generalized}, we have 
\begin{lemma}\label{lma:Lemma 12 Gen Log Hiai}
Let $\mathbf{x}, \mathbf{y} \in \mathbb{R}^r_{\geq 0}$ such that  $x_1 \geq \cdots \geq x_r$ and $y_1 \geq \cdots \geq y_r$ with $\mathbf{x}\prec_{\log} \mathbf{y}$. Also let $\mathbf{y}_i = [y_{i;1}, \cdots , y_{i;r} ] \in \mathbb{R}^r_{\geq 0}$ be a sequence of vectors such that $y_{i;1} \geq \cdots \geq y_{i;r} > 0$ and $\mathbf{y}_i \rightarrow \mathbf{y}$ as $i \rightarrow \infty$. Then, there exists $i_0 \in \mathbb{N}$ and $\mathbf{x}_i  = [x_{i;1}, \cdots , x_{i;r} ] \in \mathbb{R}^r_{\geq 0}$ for $i \geq i_0$ such that $x_{i;1} \geq \cdots \geq x_{i;r} > 0$, $\mathbf{x}_i \rightarrow \mathbf{x}$ as $i \rightarrow \infty$, and 
\begin{eqnarray}
\mathbf{x}_i \prec_{\log} \mathbf{y}_i \mbox{~~for $i \geq i_0$.} 
\end{eqnarray}
\end{lemma}

For any function $f$ on $\mathbb{R}_{\geq 0}$, the term $f(\mathbf{x}$ is defined as $f(\mathbf{x}) \define (f(x_1), \cdots, f(x_r))$ with conventions $e^{ - \infty} = 0$ and $\log 0 = - \infty$. 

\section{Multivariate Tensor Norm Inequalities}\label{sec:Multivariate Tensor Norm Inequalities}

In this section, we will develop several theorems about majorization in Section~\ref{sec:Majorization wtih Integral Average}, and log majorization wtih integral average in Section~\ref{sec:Log-Majorization wtih Integral Average}. These majorization related theorems will provide tools to establish new inequalities for unitarily norms of multivariate tensors in Section~\ref{sec:Multivaraite Tensor Norm Inequalities}.

\subsection{Majorization wtih Integral Average}\label{sec:Majorization wtih Integral Average}


Let $\Omega$ be a $\sigma$-compact metric space and $\nu$ a probability measure on the Boreal $\sigma$-field of $\Omega$. Let $\mathcal{C}, \mathcal{D}_\tau \in \mathbb{C}^{I_1 \times \cdots \times I_N \times I_1 \times \cdots \times I_N}$ be Hermitian tensors with Hermitian rank $r$. We further assume that tensors $\mathcal{C}, \mathcal{D}_\tau$ are uniformly bounded in their norm for $\tau \in \Omega$. Let $\tau \in\Omega \rightarrow  \mathcal{D}_\tau$ be a continuous function such that $\sup \{\left\Vert  D_{\tau} \right\Vert: \tau \in \Omega  \} < \infty$. For notational convenience, we define the following relation:
\begin{eqnarray}\label{eq:integral eigen vector rep}
\left[ \int_{\Omega} \lambda_1(\mathcal{D}_\tau) d\nu(\tau), \cdots, \int_{\Omega} \lambda_r(\mathcal{D}_\tau) d\nu(\tau) \right] \define \int_{\Omega^r} \vec{\lambda}(\mathcal{D}_\tau) d\nu^r(\tau).
\end{eqnarray}
If $f$ is a single variable function, the notation $f(\mathcal{C})$ represents a tensor function with respect to the tensor $\mathcal{C}$. 

\begin{theorem}\label{thm:weak int average thm 4}
Let $\Omega, \nu, \mathcal{C}, \mathcal{D}_\tau$ be defined as the beginning part of Section~\ref{sec:Majorization wtih Integral Average}, and $f: \mathbb{R} \rightarrow [0, \infty)$ be a non-decreasing convex function, we have following two equivalent statements:
\begin{eqnarray}\label{eq1:thm:weak int average thm 4}
\vec{\lambda}(\mathcal{C}) \prec_w  \int_{\Omega^r} \vec{\lambda}(\mathcal{D}_\tau) d\nu^r(\tau) \Longleftrightarrow \left\Vert f(\mathcal{C}) \right\Vert_{\rho} \leq 
\int_{\Omega} \left\Vert f(\mathcal{D}_{\tau}) \right\Vert_{\rho}  d\nu(\tau),
\end{eqnarray}
where $\left\Vert \cdot \right\Vert_{\rho}$ is the unitarily invariant norm defined in Eq.~\eqref{eq:def gauge func and general unitarily invariant norm}. 
\end{theorem}
\textbf{Proof:}
We assume that the left statement of Eq.~\eqref{eq1:thm:weak int average thm 4} is true and the function $f$ is a non-decreasing convex function. From Lemma~\ref{lma:Lemma 1 Gen Log Hiai}, we have 
\begin{eqnarray}\label{eq2:thm:weak int average thm 4}
\vec{\lambda}(f (\mathcal{C})) = f (\vec{\lambda}(\mathcal{C})) \prec_w  f \left(\int_{\Omega^r} \vec{\lambda}(\mathcal{D}_\tau) d\nu^r(\tau) \right).
\end{eqnarray}
From the convexity of $f$, we also have 
\begin{eqnarray}\label{eq3:thm:weak int average thm 4}
f \left(\int_{\Omega^r} \vec{\lambda}(\mathcal{D}_\tau) d\nu^r(\tau) \right) \leq \int_{\Omega^r} f(\vec{\lambda}(\mathcal{D}_\tau)) d\nu^r(\tau) = \int_{\Omega^r} \vec{\lambda} ( f(\mathcal{D}_\tau)) d\nu^r(\tau).
\end{eqnarray}
Then, we obtain $\vec{\lambda}(f (\mathcal{C}))  \prec_{w} = \int_{\Omega^r} \vec{\lambda} ( f(\mathcal{D}_\tau)) d\nu^r(\tau)$. By applying Lemma 4.4.2 in~\cite{hiai2010matrix} to both sides of $\vec{\lambda}(f (\mathcal{C}))  \prec_{w} = \int_{\Omega^r} \vec{\lambda} ( f(\mathcal{D}_\tau)) d\nu^r(\tau)$ with gauge function $\rho$, we obtain 
\begin{eqnarray}\label{eq4:thm:weak int average thm 4}
\left \Vert f(\mathcal{C}) \right\Vert_{\rho} &\leq &\rho \left( \int_{\Omega^r} \vec{\lambda} ( f(\mathcal{D}_\tau)) d\nu^r(\tau)  \right)  \nonumber \\
&\leq & \int_{\Omega} \rho(\vec{\lambda} ( f(\mathcal{D}_\tau))) d\nu(\tau) 
= \int_{\Omega} \left\Vert f(\mathcal{D}_\tau) \right\Vert_{\rho} d\nu(\tau).
\end{eqnarray}
Therefore, the right statement of Eq.~\eqref{eq1:thm:weak int average thm 4} is true from the left statement. 

On the other hand, if the right statement of Eq.~\eqref{eq1:thm:weak int average thm 4} is true, we select a function $f \define \max\{x + c, 0\} $, where $c$ is a postive real constant satisfying $\mathcal{C} + c \mathcal{I} \geq \mathcal{O}$, $\mathcal{D}_{\tau} + c \mathcal{I} \geq \mathcal{O}$ for all $\tau \in \Omega$, and tensors $\mathcal{C} + c \mathcal{I}, \mathcal{D}_{\tau} + c \mathcal{I}$ having Hermitian rank $r$. If the Ky Fan $k$-norm at the right statement of Eq.~\eqref{eq1:thm:weak int average thm 4} is applied, we have 
\begin{eqnarray}\label{eq5:thm:weak int average thm 4}
\sum\limits_{i=1}^k (\lambda_i (\mathcal{C}) + c ) \leq  
\sum\limits_{i=1}^k \int_{\Omega} ( \lambda_i (\mathcal{D}_{\tau}) + c ) d\nu(\tau).
\end{eqnarray}
Hence, $\sum\limits_{i=1}^k \lambda_i (\mathcal{C}) \leq  
\sum\limits_{i=1}^k \int_{\Omega} \lambda_i (\mathcal{D}_{\tau}) d\nu(\tau)$, this is the left statement of Eq.~\eqref{eq1:thm:weak int average thm 4}.
$\hfill \Box$

Next theorem will provide a stronger version of Theorem~\ref{thm:weak int average thm 4} by removing weak majorization conditions. 
\begin{theorem}\label{thm:weak int average thm 5}
Let $\Omega, \nu, \mathcal{C}, \mathcal{D}_\tau$ be defined as the beginning part of Section~\ref{sec:Majorization wtih Integral Average}, and $f: \mathbb{R} \rightarrow [0, \infty)$ be a convex function, we have following two equivalent statements:
\begin{eqnarray}\label{eq1:thm:weak int average thm 5}
\vec{\lambda}(\mathcal{C}) \prec  \int_{\Omega^r} \vec{\lambda}(\mathcal{D}_\tau) d\nu^r(\tau) \Longleftrightarrow \left\Vert f(\mathcal{C}) \right\Vert_{\rho} \leq 
\int_{\Omega} \left\Vert f(\mathcal{D}_{\tau}) \right\Vert_{\rho}  d\nu(\tau),
\end{eqnarray}
where $\left\Vert \cdot \right\Vert_{\rho}$ is the unitarily invariant norm defined in Eq.~\eqref{eq:def gauge func and general unitarily invariant norm}. 
\end{theorem}
\textbf{Proof:}
We assume that the left statement of Eq.~\eqref{eq1:thm:weak int average thm 5} is true and the function $f$ is a convex function. Again, from Lemma~\ref{lma:Lemma 1 Gen Log Hiai}, we have
\begin{eqnarray}\label{eq2:thm:weak int average thm 5}
\vec{\lambda}(f(\mathcal{A})) = f(\vec{\lambda}(\mathcal{A})) \prec_{w}  f \left( \left(\int_{\Omega^r} \vec{\lambda}(\mathcal{D}_\tau) d\nu^r(\tau)  \right) \right) \leq \int_{\Omega^r} f(\vec{\lambda}(\mathcal{D}_\tau)) d\nu^r(\tau),
\end{eqnarray}
then, 
\begin{eqnarray}\label{eq3:thm:weak int average thm 5}
\left\Vert f(\mathcal{A}) \right\Vert_{\rho} &\leq & \rho\left( \int_{\Omega^r} f(\vec{\lambda}(\mathcal{D}_\tau)) d\nu^r(\tau) \right) \nonumber \\
& \leq & \int_{\Omega}\rho \left( f(\vec{\lambda}(\mathcal{D}_\tau)) \right)d\nu (\tau) = 
\int_{\Omega} \left\Vert f( \mathcal{D}_\tau) \right\Vert_{\rho} d\nu (\tau). 
\end{eqnarray}
This prove the right statement of Eq.~\eqref{eq1:thm:weak int average thm 5}. 

Now, we assume that the right statement of Eq.~\eqref{eq1:thm:weak int average thm 5} is true. From Theorem~\ref{thm:weak int average thm 4}, we already have $\vec{\lambda}(\mathcal{C}) \prec_w  \int_{\Omega^r} \vec{\lambda}(\mathcal{D}_\tau) d\nu^r(\tau)$.  It is enough to prove $\sum\limits_{i=1}^r \lambda_i(\mathcal{C}) \geq \int_{\Omega} \sum\limits_{i=1}^r \lambda_i(\mathcal{D}_{\tau}) d \nu(\tau)$. We define a function $f \define \max\{c - x, 0\} $, where $c$ is a postive real constant satisfying $\mathcal{C} \leq c \mathcal{I} $, $\mathcal{D}_{\tau} \leq  c \mathcal{I}$ for all $\tau \in \Omega$ and tensors $c \mathcal{I} - \mathcal{C}, c \mathcal{I} - \mathcal{D}_{\tau}$ having Hermitian rank $r$. If the trace norm is applied, i.e., the sum of the absolute value of all eigenvalues of a Hermitian tensor, then the right statement of Eq.~\eqref{eq1:thm:weak int average thm 5} becomes
\begin{eqnarray}\label{eq4:thm:weak int average thm 5}
\sum\limits_{i=1}^r \lambda_i \left( c\mathcal{I} - \mathcal{C}\right)  \leq \int_{\Omega} 
\sum\limits_{i=1}^r \lambda_i \left( c\mathcal{I} - \mathcal{D}_{\tau}\right) d \nu(\tau).
\end{eqnarray}
The desired inequality  $\sum\limits_{i=1}^r \lambda_i(\mathcal{C}) \geq \int_{\Omega} \sum\limits_{i=1}^r \lambda_i(\mathcal{D}_{\tau}) d \nu(\tau)$ is established. $\hfill \Box$

\subsection{Log-Majorization wtih Integral Average}\label{sec:Log-Majorization wtih Integral Average}


The purpose of this section is to consider log-majorization issues for unitarily invariant norms of Hermitian tensors. In this section, let $\mathcal{C}, \mathcal{D}_\tau \in \mathbb{C}^{I_1 \times \cdots \times I_N \times I_1 \times \cdots \times I_N}$ be nonnegative Hermitian tensors with Hermitian rank $r$, i,e, all Hermitian eigenvalues are positive, and keep other notations with the same definitions as at the beginning of the Section~\ref{sec:Majorization wtih Integral Average}. For notational convenience, we define the following relation for logarithm vector:
\begin{eqnarray}\label{eq:integral eigen log vector rep}
\left[ \int_{\Omega} \log \lambda_1(\mathcal{D}_\tau) d\nu(\tau), \cdots, \int_{\Omega} \log \lambda_r(\mathcal{D}_\tau) d\nu(\tau) \right] \define \int_{\Omega^r} \log \vec{\lambda}(\mathcal{D}_\tau) d\nu^r(\tau).
\end{eqnarray}

\begin{theorem}\label{thm:weak int average thm 7}
Let $\mathcal{C}, \mathcal{D}_\tau$ be nonnegative Hermitian tensors, $f: (0, \infty) \rightarrow [0,\infty)$ be a continous function such that the mapping $x \rightarrow \log f(e^x)$ is convex on $\mathbb{R}$, and $g: (0, \infty) \rightarrow [0,\infty)$ be a continous function such that the mapping $x \rightarrow g(e^x)$ is convex on $\mathbb{R}$ , then  we have following three equivalent statements:
\begin{eqnarray}\label{eq1:thm:weak int average thm 7}
\vec{\lambda}(\mathcal{C}) &\prec_{w \log}& \exp  \int_{\Omega^r} \log \vec{\lambda}(\mathcal{D}_\tau) d\nu^r(\tau);
\end{eqnarray}
\begin{eqnarray}\label{eq2:thm:weak int average thm 7}
\left\Vert f(\mathcal{C}) \right\Vert_{\rho} &\leq &
\exp \int_{\Omega} \log \left\Vert f(\mathcal{D}_{\tau}) \right\Vert_{\rho}  d\nu(\tau);
\end{eqnarray}
\begin{eqnarray}\label{eq3:thm:weak int average thm 7}
\left\Vert g(\mathcal{C}) \right\Vert_{\rho} &\leq &
\int_{\Omega} \left\Vert g(\mathcal{D}_{\tau}) \right\Vert_{\rho}  d\nu(\tau).
\end{eqnarray}
\end{theorem}
\textbf{Proof:}
The roadmap of this proof is to prove equivalent statements between Eq.~\eqref{eq1:thm:weak int average thm 7} and Eq.~\eqref{eq2:thm:weak int average thm 7} first, followed by equivalent statements between Eq.~\eqref{eq1:thm:weak int average thm 7} and Eq.~\eqref{eq3:thm:weak int average thm 7}. 

\textbf{Eq.~\eqref{eq1:thm:weak int average thm 7} $\Longrightarrow$ Eq.~\eqref{eq2:thm:weak int average thm 7}}

There are two cases to be discussed in this part of proof: $\mathcal{C}, \mathcal{D}_\tau$ are positive Hermitian tensors, and $\mathcal{C}, \mathcal{D}_\tau$ are nonnegative Hermitian tensors. At the beginning, we consider the case that $\mathcal{C}, \mathcal{D}_\tau$ are positive Hermitian tensors.

Since $\mathcal{D}_\tau$ are positive, we can find $\varepsilon > 0$ such that $\mathcal{D}_{\tau} \geq \varepsilon \mathcal{I}$ for all $\tau \in \Omega$. From Eq.~\eqref{eq1:thm:weak int average thm 7}, the convexity of $\log f (e^x)$ and Lemma~\ref{lma:Lemma 1 Gen Log Hiai}, we have 
\begin{eqnarray}
\vec{\lambda} \left( f ( \mathcal{C})\right)  = f \left(\exp \left( \log \vec{\lambda} (\mathcal{C}) \right) \right) &\prec _w &  f \left(\exp    \int_{\Omega^r} \vec{\lambda}(\mathcal{D}_\tau) d\nu^r(\tau)             \right) \nonumber \\
& \leq &  \exp \left( \int_{\Omega^r} \log f \left( \vec{\lambda}(\mathcal{D}_\tau) \right)  d\nu^r(\tau) \right).
\end{eqnarray}
Then, from Eq.~\eqref{eq:def gauge func and general unitarily invariant norm}, we obtain
\begin{eqnarray}\label{eq4:thm:weak int average thm 7}
\left\Vert f (\mathcal{C}) \right\Vert _{\rho}
& \leq &  \rho\left(\exp \left( \int_{\Omega^r} \log  f \left( \vec{\lambda}(\mathcal{D}_\tau) \right) d\nu^r(\tau) \right)   \right).
\end{eqnarray}


From the function $f$ properties, we can assume that $f(x) > 0$ for any $x > 0$. Then, we have 
following bounded and continous maps on $\Omega$: $\tau \rightarrow \log f (\lambda_i (\mathcal{D}_{\tau}))  $ for $i \in \{1,2,\cdots, r\}$, and $\tau \rightarrow \log \left\Vert f (\mathcal{D}_{\tau}) \right\Vert_{\rho}$. Because we have $\nu (\Omega) = 1$ and $\sigma$-compactness of $\Omega$, we have $\tau_{k}^{(n)} \in \Omega$ and $\alpha_{k}^{(n)}$ for $k \in \{1,2,\cdots, n\}$ and $n \in \mathbb{N}$ with $\sum\limits_{k=1}^{n} \alpha_{k}^{(n)} = 1$ such that 
\begin{eqnarray}\label{eq:35}
\int_{\Omega} \log f (\lambda_i ( \mathcal{D}_{\tau} ) ) d \nu (\tau) = \lim\limits_{n \rightarrow \infty} \sum\limits_{k=1}^{n} \alpha_{k}^{(n)}  \log f (\lambda_i (\mathcal{D}_{\tau_k^{(n)} }))   , \mbox{for $i \in \{1,2,\cdots r \}$};
\end{eqnarray}
and 
\begin{eqnarray}\label{eq:36}
\int_{\Omega} \log \left\Vert f (\mathcal{D}_{\tau}) \right\Vert_{\rho} d \nu (\tau) = \lim\limits_{n \rightarrow \infty} \sum\limits_{k=1}^{n} \alpha_{k}^{(n)}  \log \left\Vert f (\mathcal{D}_{\tau_k^{(n)}    }) \right\Vert_{\rho} .
\end{eqnarray}
By taking the exponential at both sides of Eq.~\eqref{eq:35} and apply the gauge function $\rho$, we have
\begin{eqnarray}\label{eq:37}
\rho \left( \exp \int_{\Omega^r} \log f (\vec{\lambda} ( \mathcal{D}_{\tau} ) ) d \nu^r (\tau)  \right)= \lim\limits_{n \rightarrow \infty} \rho\left(  \prod\limits_{k=1}^{n}  f  \left( \vec{\lambda} \left(\mathcal{D}_{\tau_k^{(n)} } \right)  \right)^{ \alpha_{k}^{(n)} }  \right).
\end{eqnarray}
Similarly, by taking the exponential at both sides of Eq.~\eqref{eq:36}, we have
\begin{eqnarray}\label{eq:38}
\exp \left( \int_{\Omega} \log \left\Vert f (\mathcal{D}_{\tau}) \right\Vert_{\rho} d \nu (\tau) \right) = \lim\limits_{n \rightarrow \infty} \prod \limits_{k=1}^{n} \left\Vert f \left( \mathcal{D}_{\tau_k^{(n)}    } \right) \right\Vert^{\alpha_{k}^{(n)}}_{\rho} .
\end{eqnarray}
From Lemma~\ref{lma:Holder inquality for gauge function}, we have 
\begin{eqnarray}\label{eq:41}
\rho \left(  \prod\limits_{k=1}^{n}  f  \left( \vec{\lambda} \left(\mathcal{D}_{\tau_k^{(n)} } \right)  \right)^{ \alpha_{k}^{(n)} }  \right) & \leq & \prod \limits_{k=1}^{n} \rho^{  \alpha_{k}^{(n)}    } \left( f \left( \vec{\lambda} \left( \mathcal{D}_{ \tau_{k}^{(n)}  }\right) 
 \right) \right) \nonumber \\
&=& \prod \limits_{k=1}^{n} \rho^{  \alpha_{k}^{(n)}    } \left( \vec{\lambda}  \left( f  \left( \mathcal{D}_{ \tau_{k}^{(n)}  }\right) 
 \right) \right) \nonumber \\
&=& \prod \limits_{k=1}^{n} \left\Vert f  \left( \mathcal{D}_{ \tau_{k}^{(n)}  } \right) \right\Vert_{\rho}^{\alpha_{k}^{(n)}}
\end{eqnarray}

From Eqs.~\eqref{eq:37},~\eqref{eq:38} and~\eqref{eq:41}, we have 
\begin{eqnarray}\label{eq:42}
\rho \left( \exp \int_{\Omega^r} \log f (\vec{\lambda} ( \mathcal{D}_{\tau} ) ) d \nu^r (\tau)  \right) \leq \exp \int_{\Omega} \log \left\Vert f(\mathcal{D}_{\tau})\right\Vert_{\rho} d \nu(\tau).
\end{eqnarray}
Then, Eq.~\eqref{eq2:thm:weak int average thm 7} is proved from Eqs.~\eqref{eq4:thm:weak int average thm 7} and~\eqref{eq:42}.

Next, we consider that $\mathcal{C}, \mathcal{D}_\tau$ are nonnegative Hermitian tensors. For any $\delta > 0$, we have following log-majorization relation:
\begin{eqnarray}
\prod\limits_{i=1}^k \left( \lambda_i (\mathcal{C}) + \epsilon_{\delta} \right) 
&\leq& \prod\limits_{i=1}^k \exp  \int_{\Omega} \log \left( \lambda_i(\mathcal{D}_\tau) + \delta\right) d \nu (\tau),
\end{eqnarray}
where $\epsilon_{\delta} > 0$ and $k \in \{1,2,\cdots r \}$. Then, we can apply the previous case result about positive Hermitian tensors to positive Hermitian tensors $\mathcal{C} + \epsilon_{\delta} \mathcal{I}$ and $\mathcal{D}_\tau + \delta \mathcal{I}$, and get 
\begin{eqnarray}\label{eq:46}
\left\Vert f (\mathcal{C}) + \epsilon_{\delta}  \mathcal{I} \right\Vert_{\rho} 
&\leq& \exp \int_{\Omega} \log \left\Vert f (\mathcal{D}_{\tau}) + \delta  \mathcal{I} \right\Vert_{\rho} 
d \nu (\tau)
\end{eqnarray}
As $\delta \rightarrow 0$, Eq.~\eqref{eq:46} will give us Eq.~\eqref{eq2:thm:weak int average thm 7} for nonnegative Hermitian tensors.  

\textbf{Eq.~\eqref{eq1:thm:weak int average thm 7} $\Longleftarrow$ Eq.~\eqref{eq2:thm:weak int average thm 7}}

We consider positive Hermitian tensors at first phase by assuming that $\mathcal{D}_{\tau}$ are 
positive Hermitian for all $\tau \in \Omega$. We may also assume that the tensor $\mathcal{C}$ is a positive Hermitian tensor. Since if this is a nonnegative Hermitian tensor, i.e., some $\lambda_i = 0$, we always have following inequality valid:
\begin{eqnarray}
\prod\limits_{i=1}^k \lambda_i (\mathcal{C}) \leq \prod\limits_{i=1}^k 
\exp \int_{\Omega} \log \lambda_i (\mathcal{D}_{\tau}) d \nu (\tau)
\end{eqnarray}

If we apply $f(x) = x^p$ for $p > 0$ and $\left\Vert \cdot \right\Vert_{\rho}$ as Ky Fan $k$-norm in Eq.~\eqref{eq2:thm:weak int average thm 7}, we have 
\begin{eqnarray}\label{eq:50}
\log \sum\limits_{i=1}^k \lambda^p_i \left(\mathcal{C}\right) \leq \int_{\Omega} \log \sum\limits_{i=1}^k \lambda_i^p\left( \mathcal{D}_{\tau} \right) d \nu (\tau).
\end{eqnarray}
If we add $\log \frac{1}{k}$ and multiply $\frac{1}{p}$ at both sides of Eq.~\eqref{eq:50}, we have 
\begin{eqnarray}\label{eq:51}
\frac{1}{p}\log \left( \frac{1}{k}\sum\limits_{i=1}^k \lambda^p_i \left(\mathcal{C}\right) \right)\leq \int_{\Omega} \frac{1}{p} \log \left( \frac{1}{k}\sum\limits_{i=1}^k \lambda_i^p\left( \mathcal{D}_{\tau} \right) \right) d \nu (\tau).
\end{eqnarray}
From L'Hopital's Rule, if $p \rightarrow 0$, we have 
\begin{eqnarray}\label{eq:52}
\frac{1}{p}\log \left( \frac{1}{k}\sum\limits_{i=1}^k \lambda^p_i \left(\mathcal{C}\right) \right) \rightarrow \frac{1}{k} \sum\limits_{i=1}^k \log \lambda_i (\mathcal{C}),
\end{eqnarray}
and 
\begin{eqnarray}\label{eq:53}
\frac{1}{p}\log \left( \frac{1}{k}\sum\limits_{i=1}^k \lambda^p_i \left(\mathcal{D}_{\tau}\right) \right) \rightarrow \frac{1}{k} \sum\limits_{i=1}^k \log \lambda_i (\mathcal{D}_{\tau}),
\end{eqnarray}
where $\tau \in \Omega$. Appling Eqs.~\eqref{eq:52} and~\eqref{eq:53} into Eq.~\eqref{eq:51} and taking $p \rightarrow 0$, we have 
\begin{eqnarray}
\sum\limits_{i=1}^k \lambda_i (\mathcal{C}) \leq \int_{\Omega} \sum\limits_{i=1}^k 
 \log \lambda_i (\mathcal{D}_{\tau}) d \nu (\tau).
\end{eqnarray}
Therefore, Eq.~\eqref{eq1:thm:weak int average thm 7} is true for positive Hermitian tensors. 

For nonnegative Hermitian tensors $\mathcal{D}_{\tau}$, since Eq.~\eqref{eq2:thm:weak int average thm 7} is valid for $\mathcal{D}_{\tau} + \delta \mathcal{I}$ for any $\delta > 0$, we can apply the previous case result about positive Hermitian tensors to $\mathcal{D}_{\tau} + \delta \mathcal{I}$ and obtain
\begin{eqnarray}
\prod\limits_{i=1}^k \lambda_i (\mathcal{C}) \leq \prod\limits_{i=1}^k  \exp 
\int_{\Omega}   \log  \left( \lambda_i (\mathcal{D}_{\tau}) + \delta \right) d \nu (\tau),
\end{eqnarray}
where $k \in \{1,2,\cdots, r\}$. Eq.~\eqref{eq1:thm:weak int average thm 7} is still true for nonnegative Hermitian tensors as $\delta \rightarrow 0$.

\textbf{Eq.~\eqref{eq1:thm:weak int average thm 7} $\Longrightarrow$ Eq.~\eqref{eq3:thm:weak int average thm 7}}

If $\mathcal{C}, \mathcal{D}_\tau$ are positive Hermitian tensors, and $\mathcal{D}_{\tau} \geq \delta \mathcal{I}$ for all $\tau \in \Omega$. From Eq.~\eqref{eq1:thm:weak int average thm 7}, we have 
\begin{eqnarray}
\vec{\lambda} (\log \mathcal{C}) = \log \vec{\lambda}(\mathcal{C}) \prec_{w}
\int_{\Omega^r} \log \vec{\lambda}(\mathcal{D}_{\tau}) d \nu^r (\tau) = 
\int_{\Omega^r} \vec{\lambda}( \log \mathcal{D}_{\tau}) d \nu^r (\tau).
\end{eqnarray}
If we apply Theorem~\ref{thm:weak int average thm 4} to $\log \mathcal{C}$, $\log \mathcal{D}_{\tau}$ with function $f(x) = g(e^x)$, where $g$ is used in Eq.~\eqref{eq3:thm:weak int average thm 7}, Eq.~\eqref{eq3:thm:weak int average thm 7} is implied. 

If $\mathcal{C}, \mathcal{D}_\tau$ are nonnegative Hermitian tensors and any $\delta > 0$, we can find $\epsilon_{\delta} \in (0, \delta)$ to satisfy following:
\begin{eqnarray}\label{eq:45}
\prod\limits_{i=1}^k\left(\lambda_i(\mathcal{C}) + \epsilon_{\delta}\right) \leq 
\prod\limits_{i=1}^k \exp \int_{\Omega}   \log \left( \lambda_i(\mathcal{D}_{\tau}) + \delta  \right) d \nu (\tau).
\end{eqnarray}
Then, from positive Hermitian tensor case, we have 
\begin{eqnarray}\label{eq:45-1}
\left\Vert g( \mathcal{C} + \epsilon_{\delta} \mathcal{I} ) \right\Vert_{\rho}
\leq \int_{\Omega} \left\Vert   g( \mathcal{D}_{\tau} + \delta \mathcal{I} )    \right\Vert_{\rho}
d \nu (\tau).
\end{eqnarray}
Eq.~\eqref{eq3:thm:weak int average thm 7} is obtained by taking $\delta \rightarrow 0$ in Eq.~\eqref{eq:45-1}. 

\textbf{Eq.~\eqref{eq1:thm:weak int average thm 7} $\Longleftarrow$ Eq.~\eqref{eq3:thm:weak int average thm 7}}

For $k \in \{1,2,\cdots, r \}$, if we apply $g(x) = \log (\delta + x )$, where $\delta >0$, and Ky Fan $k$-norm in Eq.~\eqref{eq3:thm:weak int average thm 7}, we have 
\begin{eqnarray}
\sum\limits_{i=1}^k \log \left(\delta + \lambda_i \left(\mathcal{C} \right) \right)
\leq \sum\limits_{i=1}^k \int_{\Omega} \log \left( \delta + \lambda_{i}(\mathcal{D}_{\tau}) \right) d \nu (\tau).
\end{eqnarray}
Then, we have following relation as $\delta \rightarrow 0$:
\begin{eqnarray}
\sum\limits_{i=1}^k \log \lambda_i \left(\mathcal{C} \right) 
\leq \sum\limits_{i=1}^k \int_{\Omega} \log  \lambda_{i}(\mathcal{D}_{\tau}) d \nu (\tau).
\end{eqnarray}
Therefore, Eq.~\eqref{eq1:thm:weak int average thm 7} ccan be derived from Eq.~\eqref{eq3:thm:weak int average thm 7}. $\hfill \Box$

Next theorem will extend Theorem~\ref{thm:weak int average thm 7} to non-weak version.

\begin{theorem}\label{thm:int log average thm 10}
Let $\mathcal{C}, \mathcal{D}_\tau$ be nonnegative Hermitian tensors with $\int_{\Omega} \left\Vert \mathcal{D}_{\tau}^{-p}\right\Vert_\rho d \nu (\tau) < \infty$ for any $p > 0$, $f: (0, \infty) \rightarrow [0,\infty)$ be a continous function such that the mapping $x \rightarrow \log f(e^x)$ is convex on $\mathbb{R}$, and $g: (0, \infty) \rightarrow [0,\infty)$ be a continous function such that the mapping $x \rightarrow g(e^x)$ is convex on $\mathbb{R}$ , then  we have following three equivalent statements:
\begin{eqnarray}\label{eq1:thm:int average thm 10}
\vec{\lambda}(\mathcal{C}) &\prec_{\log}& \exp  \int_{\Omega^r} \log \vec{\lambda}(\mathcal{D}_\tau) d\nu^r(\tau);
\end{eqnarray}
\begin{eqnarray}\label{eq2:thm:int average thm 10}
\left\Vert f(\mathcal{C}) \right\Vert_{\rho} &\leq &
\exp \int_{\Omega} \log \left\Vert f(\mathcal{D}_{\tau}) \right\Vert_{\rho}  d\nu(\tau);
\end{eqnarray}
\begin{eqnarray}\label{eq3:thm:int average thm 10}
\left\Vert g(\mathcal{C}) \right\Vert_{\rho} &\leq &
\int_{\Omega} \left\Vert g(\mathcal{D}_{\tau}) \right\Vert_{\rho}  d\nu(\tau).
\end{eqnarray}
\end{theorem}
\textbf{Proof:}

The proof plan is similar to the proof in Theorem~\ref{thm:weak int average thm 7}. We prove the equivalence between Eq.~\eqref{eq1:thm:int average thm 10} and Eq.~\eqref{eq2:thm:int average thm 10} first, then prove the equivalence between Eq.~\eqref{eq1:thm:int average thm 10} and Eq.~\eqref{eq3:thm:int average thm 10}. 

\textbf{Eq.~\eqref{eq1:thm:int average thm 10} $\Longrightarrow$ Eq.~\eqref{eq2:thm:int average thm 10}}

First, we assume that $\mathcal{C}, \mathcal{D}_\tau$ are postiive Hermitian tensors with $\mathcal{D}_\tau \geq \delta \mathcal{I}$ for all $\tau \in \Omega$. The corresponding part of the proof in Theorem~\ref{thm:weak int average thm 7} about positive Hermitian tensors $\mathcal{C}, \mathcal{D}_\tau$ can be applied here. 

For case that $\mathcal{C}, \mathcal{D}_\tau$ are nonnegative Hermitian tensors, we have 
\begin{eqnarray}
\prod\limits_{i=1}^k \lambda_i (\mathcal{C}) \leq \prod\limits_{i=1}^k 
\exp \int_{\Omega} \log \left( \lambda_i (\mathcal{D}_{\tau}) + \delta_n \right)d \nu (\tau), 
\end{eqnarray}
where $\delta_n > 0$ and $\delta_n \rightarrow 0$. Because $\int_{\Omega^r} \log \left( \vec{\lambda} (\mathcal{D}_{\tau}) + \delta_n \right) d \nu^r (\tau) \rightarrow \int_{\Omega^r} \log \vec{\lambda} (\mathcal{D}_\tau)  d \nu^r (\tau) $ as $n \rightarrow \infty$, from Lemma~\ref{lma:Lemma 12 Gen Log Hiai}, we can find $\mathbf{a}^{(n)}$ with $n \geq n_0$ such that $a^{(n)}_1 \geq \cdots \geq a^{(n)}_r > 0$, $\mathbf{a}^{(n)} \rightarrow \vec{\lambda}(\mathcal{C})$ and $ \mathbf{a}^{(n)} \prec_{\log}  \exp \int_{\Omega^r} \log \vec{\lambda} \left(\mathcal{D}_{\tau} + \delta_n \mathcal{I} \right) d \nu^r (\tau)$

Selecting $\mathcal{C}^{(n)}$ with $\vec{\lambda} ( \mathcal{C}^{(n)})  = \mathbf{a}^{(n)} $ and applying positive Hermitian tensors case to $\mathcal{C}^{(n)}$ and $\mathcal{D}_{\tau} + \delta_n \mathcal{I}$, we obtain
\begin{eqnarray}\label{eq:74}
\left\Vert f (\mathcal{C}^{(n)}) \right\Vert_{\rho} \leq \exp \int_{\Omega} \log \left\Vert f (\mathcal{D}_{\tau} + \delta_n \mathcal{I}) \right\Vert_{\rho} d \nu (\tau)
\end{eqnarray}
where $n \geq n_0$.

There are two situations for the function $f$ near $0$: $f(0^{+}) < \infty$ and $f(0^{+}) = \infty$. For the case with $f(0^{+}) < \infty$, we have 
\begin{eqnarray}\label{eq:75-1}
\left\Vert f (\mathcal{C}^{(n)} )\right\Vert_{\rho} = \rho( f (\mathbf{a}^{(n)}))
\rightarrow \rho (f ( \vec{\lambda}(\mathcal{C}))) = \left\Vert f (\mathcal{C})\right\Vert_{\rho}, 
\end{eqnarray}
and
\begin{eqnarray}\label{eq:75-2}
\left\Vert f (\mathcal{D}_{\tau} + \delta_n \mathcal{I} )\right\Vert_{\rho} 
\rightarrow \left\Vert f (\mathcal{D}_{\tau})\right\Vert_{\rho}, 
\end{eqnarray}
where $\tau \in \Omega$ and $n \rightarrow \infty$. From Fatou–Lebesgue theorem, we then have 
\begin{eqnarray}\label{eq:76}
\limsup\limits_{n \rightarrow \infty} \int_{\Omega} \log \left\Vert f (\mathcal{D}_{\tau} + \delta_n \mathcal{I} )\right\Vert_{\rho} d \nu (\tau) \leq \int_{\Omega} \log \left\Vert f(\mathcal{D}_{\tau}) \right\Vert_{\rho}.
\end{eqnarray}
By taking $n \rightarrow \infty$ in Eq.~\eqref{eq:74} and using Eqs.~\eqref{eq:75-1},~\eqref{eq:75-2},~\eqref{eq:76}, we have Eq.~\eqref{eq2:thm:int average thm 10} for case that $f(0^{+}) < \infty$.

For the case with $f(0^{+}) = \infty$, we assume that $\int_{\Omega} \log \left\Vert f (\mathcal{D}_{\tau}) \right\Vert_{\rho} d \nu (\tau) < \infty$ (since the inequality in Eq.~\eqref{eq2:thm:int average thm 10} is always true for $\int_{\Omega} \log \left\Vert f (\mathcal{D}_{\tau}) \right\Vert_{\rho} d \nu (\tau) = \infty$). Since $f$ is decreasing on $(0, \epsilon)$ for some $\epsilon > 0$. We claim that the following relation is valid: there are two constants $a, b > 0$ such that 
\begin{eqnarray}\label{eq:77}
a \leq \left\Vert f (\mathcal{D}_{\tau} + \delta_n \mathcal{I}) \right\Vert_{\rho}
\leq \left\Vert f (\mathcal{D}_{\tau}) \right\Vert_{\rho} + b,
\end{eqnarray}
for all $\tau \in \Omega$ and $n \geq n_0$. If Eq.~\eqref{eq:77} is valid and $\int_{\Omega} \log \left\Vert f (\mathcal{D}_{\tau}) \right\Vert_{\rho} d \nu (\tau) < \infty$, from Lebesgue's dominated convergence theorem, we also have Eq.~\eqref{eq2:thm:int average thm 10} for case that $f(0^{+}) = \infty$ by taking $n \rightarrow \infty$ in Eq.~\eqref{eq:74}. 

Below, we will prove the claim stated by Eq.~\eqref{eq:77}. By the uniform boundedness of tensors $\mathcal{D}_{\tau}$, there is a constant $\kappa >0$ such that 
\begin{eqnarray}
0 < \mathcal{D}_{\tau} + \delta_n \mathcal{I} \leq \kappa \mathcal{I},
\end{eqnarray}
where $\tau \in \Omega$ and $ n \geq n_0$. We may assume that $\mathcal{D}_\tau$ is positive Hermitian tensors because $\left\Vert f (\mathcal{D}_{\tau}) \right\Vert_{\rho} = \infty$, i.e., Eq.~\eqref{eq:77} being true automatically, when $\mathcal{D}_\tau$ is nonnegative Hermitian tensors. From Eq.~\eqref{eq:Hermitian Eigen Decom}, we have 
\begin{eqnarray}
f(\mathcal{D}_{\tau} + \delta_n \mathcal{I}) &=&  \sum\limits_{i', \mbox{s.t. $\lambda_{i'}(\mathcal{D}_{\tau})+  \delta_n < \epsilon$}} f(\lambda_{i'}(\mathcal{D}_{\tau}) + \delta_n ) \mathcal{U}_{i'} \otimes \mathcal{U}^{H}_{i'} + \nonumber \\
&  & \sum\limits_{j', \mbox{s.t. $\lambda_{j'}(\mathcal{D}_{\tau}) +  \delta_n 
\geq \epsilon$}} f(\lambda_{j'}(\mathcal{D}_{\tau}) + \delta_n ) \mathcal{U}_{j'} \otimes \mathcal{U}^{H}_{j'}  \nonumber \\
&\leq &  \sum\limits_{i', \mbox{s.t. $\lambda_{i'}(\mathcal{D}_{\tau}) +  \delta_n < \epsilon$}} f(\lambda_{i'}(\mathcal{D}_{\tau}) ) \mathcal{U}_{i'} \otimes \mathcal{U}^{H}_{i'} + \nonumber \\
&  & \sum\limits_{j', \mbox{s.t. $\lambda_{j'}(\mathcal{D}_{\tau})  +  \delta_n 
\geq \epsilon$}} f(\lambda_{j'}(\mathcal{D}_{\tau}) + \delta_n ) \mathcal{U}_{j'} \otimes \mathcal{U}^{H}_{j'}  \nonumber \\
&\leq & f(\mathcal{D}_{\tau}) + \sum\limits_{j', \mbox{s.t. $\lambda_{j'}(\mathcal{D}_{\tau}) +  \delta_n \geq \epsilon$}} f(\lambda_{j'}(\mathcal{D}_{\tau}) + \delta_n ) \mathcal{U}_{j'} \otimes \mathcal{U}^{H}_{j'}.
\end{eqnarray}
Therefore, the claim in Eq.~\eqref{eq:77} follows by the triangle inequality for $\left\Vert \cdot \right\Vert_{\rho}$ and $f(\lambda_{j'}(\mathcal{D}_{\tau}) + \delta_n )  < \infty$ for $\lambda_{j'}(\mathcal{D}_{\tau}) +  \delta_n \geq \epsilon$. 

\textbf{Eq.~\eqref{eq1:thm:int average thm 10} $\Longleftarrow$ Eq.~\eqref{eq2:thm:int average thm 10}}

The weak majorization relation 
\begin{eqnarray}\label{eq:82}
\prod\limits_{i=1}^{k} \lambda_i (\mathcal{C}) \leq \prod\limits_{i=1}^{k} \exp \int_{\Omega} \log \lambda_i (\mathcal{D}_\tau) d \nu (\tau),
\end{eqnarray}
is valid for $k < r$ from Eq.~\eqref{eq1:thm:weak int average thm 7} $\Longrightarrow$ Eq.~\eqref{eq2:thm:weak int average thm 7} in Theorem~\ref{thm:weak int average thm 7}.   We wish to prove that Eq.~\eqref{eq:82} becomes equal for $k = r$. It is equivalent to prove that 
\begin{eqnarray}\label{eq:83}
\log \det\nolimits_H (\mathcal{C}) \geq   \int_{\Omega} \log \det\nolimits_H (\mathcal{D}_{\tau}) d \nu (\tau),
\end{eqnarray}
where $\det\nolimits_H ( \cdot )$ is the \emph{Hermitian determinant}. We can assume that $ \int_{\Omega} \log \det\nolimits_H (\mathcal{D}_{\tau}) d \nu (\tau) \geq - \infty$ since Eq.~\eqref{eq:83} is true for  $ \int_{\Omega} \log \det\nolimits_H (\mathcal{D}_{\tau}) d \nu (\tau) = - \infty$. Then, $\mathcal{D}_\tau$ are positive Hermitian tensors. 

If we scale tensors $\mathcal{C}, \mathcal{D}_{\tau}$ as $a \mathcal{C}, a\mathcal{D}_{\tau}$ by some $a >0$, we can assume $\mathcal{D}_{\tau} \leq \mathcal{I}$ and $\lambda_i(\mathcal{D}_\tau) \leq 1$ for all $ \tau \in \Omega$ and $ i \in \{1,2,\cdots, r\}$. Then for any $p >0$, we have 
\begin{eqnarray}
\frac{1}{r} \left\Vert \mathcal{D}_{\tau}^{-p} \right\Vert_1 \leq \lambda^{-p}_r (\mathcal{D}_{\tau} ) \leq ( \det\nolimits_H  (\mathcal{D}_{\tau})  )^{-p},
\end{eqnarray}
and 
\begin{eqnarray}\label{eq:85}
\frac{1}{p} \log \left( \frac{\left\Vert \mathcal{D}^{-p}_\tau \right\Vert_1 }{r}\right)
\leq - \log \det\nolimits_H  (\mathcal{D}_{\tau}). 
\end{eqnarray}
If we use tensor trace norm, represented by $\left\Vert \cdot \right\Vert_1$,  as unitarily invariant tensor norm and $f(x) = x^{-p}$ for any $p > 0$ in Eq.~\eqref{eq2:thm:int average thm 10}, we obtain
\begin{eqnarray}\label{eq:86}
\log \left\Vert \mathcal{C}^{-p} \right\Vert_1 \leq \int_{\Omega} \log \left\Vert \mathcal{D}^{-p}_{\tau} \right\Vert_1 d \nu(\tau).
\end{eqnarray}
By adding $\log \frac{1}{r}$ and multiplying $\frac{1}{p}$ for both sides of Eq.~\eqref{eq:86}, we have 
\begin{eqnarray}\label{eq:87}
\frac{1}{p}\log \left( \frac{\left\Vert \mathcal{C}^{-p} \right\Vert_1 }{r} \right)
\leq \int_{\Omega}\frac{1}{p} \log \left(  \frac{\left\Vert \mathcal{D}_\tau^{-p} \right\Vert_1 }{r}           \right) d \nu (\tau)
\end{eqnarray}
Similar to Eqs.~\eqref{eq:52} and~\eqref{eq:53}, we have following two relations as $p \rightarrow 0$:
\begin{eqnarray}\label{eq:88}
\frac{1}{p}\log \left( \frac{\left\Vert \mathcal{C}^{-p} \right\Vert_1 }{r} \right) \rightarrow \frac{- 1}{r} \log \det\nolimits_H (\mathcal{C}),
\end{eqnarray}
and 
\begin{eqnarray}\label{eq:89}
\frac{1}{p}\log \left( \frac{\left\Vert \mathcal{D}_{\tau}^{-p} \right\Vert_1}{r}  \right) \rightarrow \frac{- 1}{r} \log \det\nolimits_H (\mathcal{D}_\tau).
\end{eqnarray}
From Eq.~\eqref{eq:85} and Lebesgue's dominated convergence theorem, we have 
\begin{eqnarray}\label{eq:90}
\lim\limits_{p  \rightarrow 0} \int_{\Omega}\frac{1}{p}\log \left( \frac{\left\Vert \mathcal{D}_{\tau}^{-p} \right\Vert_1}{r}  \right) d \nu (\tau)= \frac{-1}{r} \int_{\Omega} \log \det\nolimits_{H}(\mathcal{D}_{\tau})     \nu (\tau) 
\end{eqnarray}
Finally, we have Eq.~\eqref{eq:83} from Eqs.~\eqref{eq:87} and~\eqref{eq:90}.

\textbf{Eq.~\eqref{eq1:thm:int average thm 10} $\Longrightarrow$ Eq.~\eqref{eq3:thm:int average thm 10}}

First, we assume that $\mathcal{C}, \mathcal{D}_\tau$ are positive Hermitian tensors and $\mathcal{D}_{\tau} \geq \delta \mathcal{I}$ for $\tau \in \Omega$. From Eq.~\eqref{eq1:thm:int average thm 10}, we can apply Theorem~\ref{thm:weak int average thm 5} to $\log \mathcal{C}, \log \mathcal{D}_\tau$ and $f(x) = g(e^x)$ to obtain Eq.~\eqref{eq3:thm:int average thm 10}. 

For $\mathcal{C}, \mathcal{D}_\tau$ are nonnegative Hermitian tensors, we can choose  $\mathbf{a}^{(n)}$ and corresponding  $\mathcal{C}^{(n)}$ for $n \geq n_0$ given $\delta_n \rightarrow 0$ with $\delta_n > 0$ as the proof in Eq.~\eqref{eq1:thm:int average thm 10} $\Longrightarrow$ Eq.~\eqref{eq2:thm:int average thm 10}. Since tensors $\mathcal{C}^{(n)}, \mathcal{D}_\tau + \delta_n \mathcal{I}$ are postive Hermitian tensors, we then have 
\begin{eqnarray}\label{eq:92}
\left \Vert g ( \mathcal{C}^{(n)}  ) \right\Vert_{\rho} \leq \int_{\Omega} \left\Vert g ( \mathcal{D}_\tau + \delta_n \mathcal{I} ) \right\Vert_{\rho} d \nu (\tau).
\end{eqnarray}
If $g(0^+) < \infty$, Eq.~\eqref{eq3:thm:int average thm 10} is obtained from Eq.~\eqref{eq:92} by taking $n \rightarrow \infty$. On the other hand, if $g(0^+) = \infty$, we can apply the argument similar to the portion about $f(0^+) = \infty$ in the proof for Eq.~\eqref{eq1:thm:int average thm 10} $\Longrightarrow$ Eq.~\eqref{eq2:thm:int average thm 10} to get $a, b > 0$ such that 
\begin{eqnarray}\label{eq:92 infty}
a \leq \left\Vert g ( \mathcal{D}_\tau + \delta_n \mathcal{I} ) \right\Vert_\rho \leq   \left\Vert  g ( \mathcal{D}_\tau ) \right\Vert_\rho + b,
\end{eqnarray}
for all $\tau \in \Omega$ and $n \geq n_0$. Since the case that $\int_{\Omega} \left\Vert  g ( \mathcal{D}_\tau ) \right\Vert_{\rho} d \nu (\tau) = \infty$ will have Eq.~\eqref{eq3:thm:int average thm 10}, we only consider the case that $\int_{\Omega} \left\Vert  g ( \mathcal{D}_\tau ) \right\Vert_{\rho} d \nu (\tau) < \infty$. Then, we have Eq.~\eqref{eq3:thm:int average thm 10} from Eqs.~\eqref{eq:92},~\eqref{eq:92 infty} and Lebesgue's dominated convergence theorem.

\textbf{Eq.~\eqref{eq1:thm:int average thm 10} $\Longleftarrow$ Eq.~\eqref{eq3:thm:int average thm 10}}

The weak majorization relation
\begin{eqnarray}\label{eq:94}
\sum\limits_{i=1}^{k} \log \lambda_i (\mathcal{C}) \leq 
\sum\limits_{i=1}^{k} \int_{\Omega} \log \lambda_i (\mathcal{D}_\tau) d \nu (\tau)
\end{eqnarray}
is true from the implication from Eq.~\eqref{eq1:thm:weak int average thm 7} to Eq.~\eqref{eq3:thm:weak int average thm 7} in Theorem~\ref{thm:weak int average thm 7}. We have to show that this relation becomes identity for $k=r$. If we apply $\left\Vert \cdot  \right\Vert_{\rho} = \left\Vert \cdot \right\Vert_1$ and $g(x) = x^{-p}$ for any $p > 0$ in Eq.~\eqref{eq3:thm:int average thm 10}, we have 
\begin{eqnarray}\label{eq:95}
\frac{1}{p} \log \left( \frac{ \left\Vert \mathcal{C}^{-p}\right\Vert_1  }{r} \right)
\leq  \frac{1}{p} \log \left( \int_{\Omega} \frac{\left\Vert \mathcal{D}_\tau^{-p}\right\Vert_1}{r}   d \nu(\tau)  \right).
\end{eqnarray}
Then, we will get 
\begin{eqnarray}\label{eq:96}
\frac{- \log \det\nolimits_H (\mathcal{C})}{r} &=& \lim\limits_{p \rightarrow 0} \frac{1}{p} \log \left( \frac{ \left\Vert \mathcal{C}^{-p}\right\Vert_1  }{r} \right) \nonumber \\
& \leq & \lim\limits_{p \rightarrow 0} \frac{1}{p} \log \left( \int_{\Omega} \frac{\left\Vert \mathcal{D}_\tau^{-p}\right\Vert_1}{r}   d \nu(\tau)  \right)=_1 \frac{ - \int_{\Omega} \log \det\nolimits_H (\mathcal{D}_{\tau}) d \nu (\tau)  }{r},
\end{eqnarray}
which will prove the identity for Eq.~\eqref{eq:94} when $k = r$. The equality in $=_1$ will be proved by the following Lemma~\ref{lma:15}.
$\hfill \Box$

\begin{lemma}\label{lma:15}
Let $\mathcal{D}_\tau$ be nonnegative Hermitian tensors with $\int_{\Omega} \left\Vert \mathcal{D}_{\tau}^{-p}\right\Vert_\rho d \nu (\tau) < \infty$ for any $p > 0$, then we have
\begin{eqnarray}\label{eq1:lma:15}
\lim\limits_{p \rightarrow 0} \left( \frac{1}{p} \log \int_{\Omega} \frac{ \left\Vert \mathcal{D}_\tau^{-p}\right\Vert_1 }{r} d \nu (\tau)\right) &=& -\frac{1}{r} \int_{\Omega} \log \det\nolimits_{H}( \mathcal{D}_{\tau} ) d \nu (\tau)
\end{eqnarray}
\end{lemma}
\textbf{Proof:}
Because $\int_{\Omega} \left\Vert \mathcal{D}_{\tau}^{-p}\right\Vert_\rho d \nu (\tau) < \infty$, we have that $\mathcal{D}_{\tau}$ are positive Hermitian tensors for $\tau$ almost everywhere in $\Omega$. Then, we have 
\begin{eqnarray}
\lim\limits_{p \rightarrow 0} \left( \frac{1}{p}\log \int_{\Omega} \frac{ \left\Vert \mathcal{D}_\tau^{-p}\right\Vert_1}{r} d \nu (\tau) \right) &=_1& \lim\limits_{p \rightarrow 0}\frac{   \int_{\Omega} \frac{ - \sum\limits_{i=1}^r  \log \lambda_i(\mathcal{D}_\tau)   }{r} d \nu (\tau)   }{  \int_{\Omega} \frac{ \left\Vert \mathcal{D}_\tau^{-p}\right\Vert_1}{r} d \nu (\tau)      } \nonumber \\
&=& \frac{-1}{r} \int_{\Omega} \sum\limits_{i=1}^r \log \lambda_i(\mathcal{D}_\tau)  d \nu (\tau)  
\nonumber \\
&=_2& \frac{-1}{r} \int_{\Omega} \log \det\nolimits_H ( \mathcal{D}_\tau ) d \nu (\tau), 
\end{eqnarray}
where $=_1$ is from L'Hopital's rule, and $=_2$ is obtained from $\det\nolimits_H$ definition.
$\hfill \Box$

\subsection{Multivaraite Tensor Norm Inequalities}\label{sec:Multivaraite Tensor Norm Inequalities}


In this section, we will apply derived majorization inequalities for tensors to multivariate tensor norm inequalities which will be used to bound random tensor concentration inequalities in later sections. We will begin to present a Lie-Trotter product formula for tensors. 
\begin{lemma}\label{lma: Lie product formula for tensors}
Let $m \in \mathbb{N}$ and $(\mathcal{L}_k)_{k=1}^{m}$ be a finite sequence of bounded tensors with dimensions $\mathcal{L}_k \in  \mathbb{C}^{I_1 \times \cdots \times I_M\times I_1 \times \cdots \times I_M}$, then we have
\begin{eqnarray}
\lim_{n \rightarrow \infty} \left(  \prod_{k=1}^{m} \exp(\frac{  \mathcal{L}_k}{n})\right)^{n}
&=& \exp \left( \sum_{k=1}^{m}  \mathcal{L}_k \right)
\end{eqnarray}
\end{lemma}
\textbf{Proof:}

We will prove the case for $m=2$, and the general value of $m$ can be obtained by mathematical induction. 
Let $\mathcal{L}_1, \mathcal{L}_2$ be bounded tensors act on some Hilbert space. Define $\mathcal{C} \define \exp( (\mathcal{L}_1 + \mathcal{L}_2)/n) $, and $\mathcal{D} \define \exp(\mathcal{L}_1/n) \star_M \exp(\mathcal{L}_2/n)$. Note we have following estimates for the norm of tensors $\mathcal{C}, \mathcal{D}$: 
\begin{eqnarray}\label{eq0: lma: Lie product formula for tensors}
\left\Vert \mathcal{C} \right\Vert, \left\Vert \mathcal{D} \right\Vert \leq \exp \left( \frac{\left\Vert \mathcal{L}_1 \right\Vert + \left\Vert \mathcal{L}_2 \right\Vert  }{n} \right) =  \left[ \exp \left(  \left\Vert \mathcal{L}_1 \right\Vert + \left\Vert \mathcal{L}_2 \right\Vert  \right) \right]^{1/n}.
\end{eqnarray}

From the Cauchy-Product formula, the tensor $\mathcal{D}$ can be expressed as:
\begin{eqnarray}
\mathcal{D} &=& \exp(\mathcal{L}_1/n) \star_M \exp(\mathcal{L}_2/n) = \sum_{i = 0}^{\infty} \frac{( \mathcal{L}_1/n)^i}{i !} \star_M \sum_{j = 0}^{\infty} \frac{( \mathcal{L}_2/n)^j}{j !} \nonumber\\
&=& \sum_{l = 0}^{\infty} n^{-l} \sum_{i=0}^l \frac{\mathcal{L}_1^i}{i!} \star_M \frac{\mathcal{L}_2^{l-i}}{(l - i)!},
\end{eqnarray}
then we can bound the norm of $\mathcal{C} - \mathcal{D}$ as 
\begin{eqnarray}\label{eq1: lma: Lie product formula for tensors}
\left\Vert \mathcal{C} - \mathcal{D} \right\Vert &=& \left\Vert \sum_{i=0}^{\infty} \frac{\left( [ \mathcal{L}_1 + \mathcal{L}_2]/n \right)^i}{i! }
 - \sum_{l = 0}^{\infty} n^{-l} \sum_{i=0}^l \frac{\mathcal{L}_1^i}{i!} \star_M \frac{\mathcal{L}_2^{l-i}}{(l - i)!} \right\Vert \nonumber \\
&=& \left\Vert \sum_{i=2}^{\infty} k^{-i} \frac{\left( [ \mathcal{L}_1 + \mathcal{L}_2] \right)^i}{i! }
 - \sum_{m = l}^{\infty} n^{-l} \sum_{i=0}^l \frac{\mathcal{L}_1^i}{i!} \star_M \frac{\mathcal{L}_2^{l-i}}{(l - i)!} \right\Vert \nonumber \\
& \leq & \frac{1}{k^2}\left[ \exp( \left\Vert \mathcal{L}_1 \right\Vert + \left\Vert \mathcal{L}_2 \right\Vert ) + \sum_{l = 2}^{\infty} n^{-l} \sum_{i=0}^l \frac{\left\Vert \mathcal{L}_1 \right\Vert^i}{i!} \cdot \frac{\left\Vert \mathcal{L}_2 \right\Vert^{l-i}}{(l - i)!} \right] \nonumber \\
& = & \frac{1}{n^2}\left[ \exp \left( \left\Vert \mathcal{L}_1 \right\Vert + \left\Vert \mathcal{L}_2 \right\Vert \right) + \sum_{l = 2}^{\infty} n^{-l} \frac{(  \left\Vert \mathcal{L}_1 \right\Vert + \left\Vert \mathcal{L}_2 \right\Vert )^l}{l!} \right] \nonumber \\
& \leq & \frac{2  \exp \left( \left\Vert \mathcal{L}_1 \right\Vert + \left\Vert \mathcal{L}_2 \right\Vert \right) }{n^2}.
\end{eqnarray}

For the difference between the higher power of $\mathcal{C}$ and $\mathcal{D}$, we can bound them as 
\begin{eqnarray}
\left\Vert \mathcal{C}^n - \mathcal{D}^n \right\Vert &=& \left\Vert \sum_{l=0}^{n-1} \mathcal{C}^m (\mathcal{C} - \mathcal{D})\mathcal{C}^{n-l-1} \right\Vert \nonumber \\
& \leq_1 &  \exp ( \left\Vert \mathcal{L}_1 \right\Vert +  \left\Vert \mathcal{L}_2 \right\Vert) \cdot n \cdot \left\Vert \mathcal{L}_1 - \mathcal{L}_2 \right\Vert,
\end{eqnarray}
where the inequality $\leq_1$ uses the following fact 
\begin{eqnarray}
\left\Vert \mathcal{C} \right\Vert^{l} \left\Vert \mathcal{D} \right\Vert^{n - l - 1} \leq \exp \left( \left\Vert \mathcal{L}_1 \right\Vert +  \left\Vert \mathcal{L}_2 \right\Vert \right)^{\frac{n-1}{n}} \leq 
 \exp\left( \left\Vert \mathcal{L}_1 \right\Vert +  \left\Vert \mathcal{L}_2 \right\Vert \right), 
\end{eqnarray}
based on Eq.~\eqref{eq0: lma: Lie product formula for tensors}. By combining with Eq.~\eqref{eq1: lma: Lie product formula for tensors}, we have the following bound
\begin{eqnarray}
\left\Vert \mathcal{C}^n - \mathcal{D}^n \right\Vert &\leq& \frac{2 \exp \left( 2  \left\Vert \mathcal{L}_1 \right\Vert  +  2  \left\Vert \mathcal{L}_2 \right\Vert \right)}{n}.
\end{eqnarray}
Then this lemma is proved when $n$ goes to infity. $\hfill \Box$

Below, new multivariate norm inequalities for tensors are provided according to previous majorization theorems. 
\begin{theorem}\label{thm:Multivaraite Tensor Norm Inequalities}
Let $\mathcal{C}_i \in \mathbb{C}^{I_1 \times \cdots \times I_N \times I_1 \times \cdots \times I_N}$ be positive Hermitian tensors for $1 \leq i \leq n$ with Hermitian rank $r$, $\left\Vert \cdot \right\Vert_{\rho}$ be a unitarily invaraint norm with corresponding gauge function $\rho$. For any continous function $f:(0, \infty) \rightarrow [0, \infty)$ such that $x \rightarrow \log f(e^x)$ is convex on $\mathbb{R}$, we have 
\begin{eqnarray}\label{eq1:thm:Multivaraite Tensor Norm Inequalities}
\left\Vert  f \left( \exp \left( \sum\limits_{i=1}^n \log \mathcal{C}_i\right)   \right)  \right\Vert_{\rho} &\leq& \exp \int_{- \infty}^{\infty} \log \left\Vert f \left( \left\vert \prod\limits_{i=1}^{n}  \mathcal{C}_i^{1 + \iota t} \right\vert\right)\right\Vert_{\rho} \beta_0(t) dt ,
\end{eqnarray}
where $\beta_0(t) = \frac{\pi}{2 (\cosh (\pi t) + 1)}$.

For any continous function $g(0, \infty) \rightarrow [0, \infty)$ such that $x \rightarrow g (e^x)$ is convex on $\mathbb{R}$, we have 
\begin{eqnarray}\label{eq2:thm:Multivaraite Tensor Norm Inequalities}
\left\Vert  g \left( \exp \left( \sum\limits_{i=1}^n \log \mathcal{C}_i\right)   \right)  \right\Vert_{\rho} &\leq& \int_{- \infty}^{\infty} \left\Vert g \left( \left\vert \prod\limits_{i=1}^{n}  \mathcal{C}_i^{1 + \iota t} \right\vert\right)\right\Vert_{\rho} \beta_0(t) dt.
\end{eqnarray}
\end{theorem}
\textbf{Proof:}
From Hirschman interpolation theorem~\cite{sutter2017multivariate} and $\theta \in [0, 1]$, we have 
\begin{eqnarray}\label{eq1:Hirschman interpolation}
\log \left\vert h(\theta) \right\vert \leq \int_{- \infty}^{\infty} \log \left\vert h(\iota t) \right\vert^{1 - \theta} \beta_{1 - \theta}(t) d t + \int_{- \infty}^{\infty} \log \left\vert h(1 +  \iota t) \right\vert^{\theta} \beta_{\theta}(t) d t , 
\end{eqnarray}
where $h(z)$ be uniformly bounded on $S \define \{ z \in \mathbb{C}: 0 \leq \Re(z) \leq 1  \}$ and holomorphic on $S$. The term $ d \beta_{\theta}(t) $ is defined as :
\begin{eqnarray}\label{eq:beta theta t def}
\beta_{\theta}(t) \define \frac{ \sin (\pi \theta)}{ 2 \theta (\cos(\pi t) + \cos (\pi \theta))  }.  
\end{eqnarray}
Let $H(z)$ be a uniformly bounded holomorphic function with values in $\mathbb{C}^{I_1 \times \cdots \times I_N \times I_1 \times \cdots \times I_N}$. Fix some $\theta \in [0, 1]$ and let $\mathcal{U}, \mathcal{V} \in \mathbb{C}^{I_1 \times \cdots \times I_N \times 1}$ be normalized tensors such that $\langle \mathcal{U}, \mathcal{H}(\theta) \star_N \mathcal{V} \rangle = \left\Vert H(\theta) \right\Vert$. If we define $h(z)$ as $h(z) \define \langle \mathcal{U}, \mathcal{H}(z) \star_N \mathcal{V} \rangle $, we have following bound: $\left\vert h(z) \right\vert \leq \left\Vert H(z) \right\Vert $ for all $z \in S$. From Hirschman interpolation theorem, we then have following interpolation theorem for tensor-valued function: 
\begin{eqnarray}\label{eq2:Hirschman interpolation}
\log \left\Vert H(\theta) \right\Vert \leq \int_{- \infty}^{\infty} \log \left\Vert H(\iota t) \right\Vert^{1 - \theta} \beta_{1 - \theta}(t) dt + \int_{- \infty}^{\infty} \log \left\Vert H(1 +  \iota t) \right\Vert^{\theta}  \beta_{\theta}(t) dt .  
\end{eqnarray}

Let $H(z) = \prod\limits_{i=1}^{n} \mathcal{C}^z_i$. Then the first term in the R.H.S. of Eq.~\eqref{eq2:Hirschman interpolation} is zero since $H(\iota t)$ is a product of unitary tensors. Then we have 
\begin{eqnarray}\label{eq:109}
\log \left\Vert \left\vert \prod\limits_{i=1}^{n} \mathcal{C}_i^{\theta} \right\vert^{\frac{1}{\theta}} \right\Vert \leq \int_{- \infty}^{\infty} \log \left\Vert    \prod\limits_{i=1}^{n} \mathcal{C}_i^{1 + \iota t}         \right\Vert \beta_{\theta}(t) d t .  
\end{eqnarray}

From Lemma~\ref{lma:antisymmetric tensor product properties}, we have following relations:
\begin{eqnarray}\label{eq:111-1}
\left\vert \prod\limits_{i=1}^n \left( \wedge^k \mathcal{C}_i \right)^{\theta}\right\vert^{\frac{1}{\theta}} = \wedge^k \left\vert  \prod\limits_{i=1}^n \mathcal{C}^{\theta}_i  \right\vert^{\frac{1}{\theta}},
\end{eqnarray}
and
\begin{eqnarray}\label{eq:111-2}
\left\vert \prod\limits_{i=1}^n \left( \wedge^k \mathcal{C}_i \right)^{1 + \iota t} \right\vert = \wedge^k \left\vert  \prod\limits_{i=1}^n \mathcal{C}^{1 + \iota t}_i  \right\vert.
\end{eqnarray}
If Eq.~\eqref{eq:109} is applied to $\wedge^k \mathcal{C}_i$ for $1 \leq k \leq r$, we have following log-majorization relation from Eqs.~\eqref{eq:111-1} and~\eqref{eq:111-2}:
\begin{eqnarray}\label{eq:112}
\log \vec{\lambda} \left(  \left\vert \prod\limits_{i=1}^{n} \mathcal{C}_i^{\theta} \right\vert^{\frac{1}{\theta}}   \right) \prec \int_{- \infty}^{\infty} \log \vec{\lambda} \left\vert \prod\limits_{i=1}^{n} \mathcal{C}_i^{1 + \iota t } \right\vert^{\frac{1}{\theta}}   \beta_{\theta}(t) d t.
\end{eqnarray}
Moreover, we have the equality condition in Eq.~\eqref{eq:112} for $k = r$ due to following identies:
\begin{eqnarray}\label{eq:113}
\det\nolimits_H \left\vert \prod\limits_{i=1}^n \mathcal{C}_i^{\theta} \right\vert^{\frac{1}{\theta}}
= \det\nolimits_H \left\vert \prod\limits_{i=1}^n \mathcal{C}_i^{1 + \iota t } \right\vert= \prod\limits_{i=1}^n \det\nolimits_H \mathcal{C}_i. 
\end{eqnarray}

At this stage, we are ready to apply Theorem~\ref{thm:int log average thm 10}  for the log-majorization provided by Eq.~\eqref{eq:112} to get following facts:
\begin{eqnarray}\label{eq:114}
\left\Vert  f \left( \left\vert \prod\limits_{i=1}^{n} \mathcal{C}_i^{\theta} \right\vert^{\frac{1}{\theta}}  \right)  \right\Vert_{\rho} &\leq& \exp \int_{- \infty}^{\infty} \log \left\Vert f \left( \left\vert \prod\limits_{i=1}^{n}  \mathcal{C}_i^{1 + \iota t} \right\vert\right)\right\Vert_{\rho} \beta_{\theta}(t) d t ,
\end{eqnarray}
and
\begin{eqnarray}\label{eq:115}
\left\Vert  g \left( \left\vert \prod\limits_{i=1}^{n} \mathcal{C}_i^{\theta} \right\vert^{\frac{1}{\theta}}  \right)  \right\Vert_{\rho} &\leq& \int_{- \infty}^{\infty} \left\Vert g \left( \left\vert \prod\limits_{i=1}^{n}  \mathcal{C}_i^{1 + \iota t} \right\vert\right)\right\Vert_{\rho}  \beta_{\theta}(t) d t.
\end{eqnarray}
From Lie product formula for tensors given by Lemma~\ref{lma: Lie product formula for tensors},  we have 
\begin{eqnarray}\label{eq:117}
 \left\vert \prod\limits_{i=1}^{n} \mathcal{C}_i^{\theta} \right\vert^{\frac{1}{\theta}}
\rightarrow \exp \left(  \sum\limits_{i=1}^{n} \log \mathcal{C}_i  \right). 
\end{eqnarray}
By setting $\theta \rightarrow 0$ in Eqs.~\eqref{eq:114},~\eqref{eq:115} and using Lie product formula given by Eq.~\eqref{eq:117},  we will get Eqs.~\eqref{eq1:thm:Multivaraite Tensor Norm Inequalities} and~\eqref{eq2:thm:Multivaraite Tensor Norm Inequalities}. 
$\hfill \Box$

\section{Generalized Tensor Chernoff and Bernstein Inequalities}\label{sec:New Tensor Inequalities}

%
%
%
%



In this section, we first utilize Theorem~\ref{thm:Multivaraite Tensor Norm Inequalities} and Laplace transform method to obtain Ky Fan $k$-norm concentration bounds for a function of tensors summation in Section~\ref{sec:Ky Fan k norm Concentration Bound}. Then, generalized tensor Chernoff and Bernstein bounds are discussed in Section~\ref{sec:Generalized Tensor Chernoff Bound} and Section~\ref{sec:Generalized Bernstein Tensor Bound}, respectively. 

\subsection{Ky Fan $k$-norm Tail Bound}\label{sec:Ky Fan k norm Concentration Bound}

We begin by introducing following two lemmas about Ky Fan $k$-norm inequalities for the product of tensors (Lemma~\ref{lma:Ky Fan Inequalities for the prod of tensosrs}) and the summation of tensors (Lemma~\ref{lma:Ky Fan Inequalities for the sum of tensosrs}).

\begin{lemma}\label{lma:Ky Fan Inequalities for the prod of tensosrs}
Let $\mathcal{C}_i \in \mathbb{C}^{I_1 \times \cdots \times I_N \times I_1 \times \cdots \times I_N}$ with Hermitian rank $r$ and let $p_i$ be positive real numbers satisfying $\sum\limits_{i=1}^m \frac{1}{p_i} = 1$. Then, we have 
\begin{eqnarray}\label{eq1:lma:Ky Fan Inequalities for the prod of tensosrs}
\left\Vert \left\vert  \prod\limits_{i=1}^{m} \mathcal{C}_i \right\vert^s \right\Vert_{(k)}
\leq  \prod\limits_{i=1}^{m} \left(  \left\Vert \left\vert \mathcal{C}_i \right\vert^{s p_i} \right\Vert_{(k)}    \right)^{\frac{1}{p_i}} 
\leq  \sum\limits_{i=1}^{m} \frac{ \left\Vert \left\vert \mathcal{C}_i \right\vert^{s p_i} \right\Vert_{(k)}      }{p_i}
\end{eqnarray}
where $s \geq 1$ and $k \in \{1,2,\cdots, r \}$. 
\end{lemma}
\textbf{Proof:}
Since we have 
\begin{eqnarray}
\left\Vert \left\vert  \prod\limits_{i=1}^{m} \mathcal{C}_i \right\vert^s \right\Vert_{(k)} 
= \sum\limits_{j=1}^k \lambda_j \left(   \left\vert  \prod\limits_{i=1}^{m} \mathcal{C}_i \right\vert^s    \right)  = \sum\limits_{j=1}^k \lambda^s_j \left(   \left\vert  \prod\limits_{i=1}^{m} \mathcal{C}_i \right\vert    \right) =
\sum\limits_{j=1}^k \sigma^s_j \left(    \prod\limits_{i=1}^{m} \mathcal{C}_i    \right),
\end{eqnarray}
where we have orders for eigenvalues as $\lambda_1 \geq \lambda_2 \geq \cdots $, and singular values as  $\sigma_1 \geq \sigma_2 \geq \cdots $. 

From Lemma~\ref{lma:antisymmetric tensor product properties}, we have 
\begin{eqnarray}
\left\Vert \left(\prod\limits_{i=1}^{m} \mathcal{C}_i \right)^{\wedge k} \right\Vert 
= \prod\limits_{j=1}^{k}  \sigma_j \left( \prod\limits_{i=1}^{m} \mathcal{C}_i  \right).
\end{eqnarray}
From the fact that the norm is submultiplicative and majorization theory, we will have 
\begin{eqnarray}\label{eq2:lma:Ky Fan Inequalities for the prod of tensosrs}
\sum\limits_{j=1}^k   \sigma^s_j \left( \prod\limits_{i=1}^{m} \mathcal{C}_i  \right)
 \leq \sum\limits_{j=1}^k \left( \prod\limits_{i=1}^{m}   \sigma^s_j (\mathcal{C}_i ) \right).
\end{eqnarray}
Then, we can apply H\"older's inequality to Eq.~\eqref{eq2:lma:Ky Fan Inequalities for the prod of tensosrs} and obtain 
\begin{eqnarray}
\sum\limits_{j=1}^k \left( \prod\limits_{i=1}^{m}   \sigma^s_j (\mathcal{C}_i ) \right) 
&\leq& \prod_{i=1}^{m} \left( \sum\limits_{j=1}^k \sigma^{s p_i}_j (\mathcal{C}_i) \right)^{\frac{1}{p_i}} = \prod_{i=1}^{m} \left( \sum\limits_{j=1}^k \lambda^{s p_i}_j ( \left\vert \mathcal{C}_i \right\vert ) \right)^{\frac{1}{p_i}}  \nonumber \\
&=& \prod_{i=1}^{m} \left( \sum\limits_{j=1}^k \lambda_j ( \left\vert \mathcal{C}_i \right\vert^{s p_i} ) \right)^{\frac{1}{p_i}}= \prod\limits_{i=1}^{m} \left(  \left\Vert \left\vert \mathcal{C}_i \right\vert^{s p_i} \right\Vert_{(k)}    \right)^{\frac{1}{p_i}}
\end{eqnarray}

The second inequality in Eq.~\eqref{eq1:lma:Ky Fan Inequalities for the prod of tensosrs} is obtained by applying Young's inequality to numbers $  \left\Vert \left\vert \mathcal{C}_i \right\vert^{s p_i} \right\Vert_{(k)} $ for $1 \leq i \leq m$. This completes the proof. $\hfill \Box$

\begin{lemma}\label{lma:Ky Fan Inequalities for the sum of tensosrs}
Let $\mathcal{C}_i \in \mathbb{C}^{I_1 \times \cdots \times I_N \times I_1 \times \cdots \times I_N}$ with Hermitian rank $r$, then we have 
\begin{eqnarray}\label{eq1:lma:Ky Fan Inequalities for the sum of tensosrs}
\left\Vert \left\vert  \sum\limits_{i=1}^{m} \mathcal{C}_i \right\vert^s \right\Vert_{(k)}
\leq  m^{s -1} \sum\limits_{i=1}^{m}  \left\Vert \left\vert \mathcal{C}_i \right\vert^{s} \right\Vert_{(k)}     
\end{eqnarray}
where $s \geq 1$ and $k \in \{1,2,\cdots, r \}$. 
\end{lemma}
\textbf{Proof:}
Since we have 
\begin{eqnarray}
\left\Vert \left\vert  \sum\limits_{i=1}^{m} \mathcal{C}_i \right\vert^s \right\Vert_{(k)}
 = \sum\limits_{j=1}^{k} \lambda_j \left( \left\vert  \sum\limits_{i=1}^{m} \mathcal{C}_i \right\vert^s \right) =  \sum\limits_{j=1}^{k} \lambda^s_j \left( \left\vert  \sum\limits_{i=1}^{m} \mathcal{C}_i \right\vert \right) = \sum\limits_{j=1}^{k} \sigma^s_j \left( \sum\limits_{i=1}^{m} \mathcal{C}_i \right). 
\end{eqnarray}
where we have orders for eigenvalues as $\lambda_1 \geq \lambda_2 \geq \cdots $, and singular values as  $\sigma_1 \geq \sigma_2 \geq \cdots $. 

From Theorem G.1.d. in~\cite{MR2759813} and Theorem 3.2 in~\cite{liang2019further}, we have Fan singular value majorization inequaliies:
\begin{eqnarray}
\sum\limits_{j = 1}^{k} \sigma_j ( \sum\limits_{i=1}^m \mathcal{C}_i ) \leq 
\sum\limits_{j = 1}^{k} \left( \sum\limits_{i=1}^m \sigma_j (\mathcal{C}_i ) \right),
\end{eqnarray}
where $k \in \{1,2,\cdots, s \}$. Then, we have 
\begin{eqnarray}
\sum\limits_{j = 1}^{k} \sigma^s_j ( \sum\limits_{i=1}^m \mathcal{C}_i ) &\leq &
\sum\limits_{j = 1}^{k} \left( \sum\limits_{i=1}^m \sigma_j (\mathcal{C}_i ) \right)^s
\leq m^{s-1} \sum\limits_{j = 1}^{k} \left( \sum\limits_{i=1}^m \sigma^s_j (\mathcal{C}_i ) \right) \nonumber \\
& = & 
m^{s-1} \sum\limits_{j = 1}^{k} \left( \sum\limits_{i=1}^m \sigma^s_j ( \left\vert \mathcal{C}_i \right\vert ) \right) = m^{s-1} \sum\limits_{j = 1}^{k} \left( \sum\limits_{i=1}^m \sigma_j ( \left\vert \mathcal{C}_i \right\vert^s ) \right) \nonumber \\
& = & m^{s -1} \sum\limits_{i=1}^{m}  \left\Vert \left\vert \mathcal{C}_i \right\vert^{s} \right\Vert_{(k)}    
\end{eqnarray}
$\hfill \Box$

Now, we are ready to present our main theorem about Ky Fan $k$-norm probability bound for a function of tensors summation. 
\begin{theorem}\label{thm:Ky Fan norm prob bound for fun of tensors sum}
Consider a sequence $\{ \mathcal{X}_j  \in \mathbb{C}^{I_1 \times \cdots \times I_N  \times I_1 \times \cdots \times I_N} \}$ of independent, random, Hermitian tensors. Let $g$ be a polynomial function with degree $n$ and nonnegative coeffecients $a_0, a_1, \cdots, a_n$ raised by power $s \geq 1$, i.e., $g(x) = \left(a_0 + a_1 x  +\cdots + a_n x^n \right)^s$. Suppose following condition is satisfied:
\begin{eqnarray}\label{eq:special cond}
g \left( \exp\left(t \sum\limits_{j=1}^{m} \mathcal{X}_j \right)\right)  \succeq \exp\left(t g \left( \sum\limits_{j=1}^{m} \mathcal{X}_j   \right) \right)~~\mbox{almost surely},
\end{eqnarray}
where $t > 0$. Then, we have 
\begin{eqnarray}\label{eq1:thm:Ky Fan norm prob bound for fun of tensors sum}
\mathrm{Pr} \left( \left\Vert g\left( \sum\limits_{j=1}^{m} \mathcal{X}_j  \right)\right\Vert_{(k)}  \geq \theta \right) \leq ~~~~~~~~~~~~~~~~~~~~~~~~~~~~~~~~~~~~~~~~~~~~ \nonumber \\
(n+1)^{s-1}\inf\limits_{t, p_j} e^{- \theta t }\left(k a_0^s + \sum\limits_{l=1}^{n} \sum\limits_{j=1}^m \frac{  a^{ls}_l  \mathbb{E} \left\Vert \exp\left( p_j  l s t \mathcal{X}_j \right) \right\Vert_{(k)} }{p_j}     \right).
\end{eqnarray}
where $\sum\limits_{j=1}^m \frac{1}{p_j} =1$ and $p_j > 0$. 
\end{theorem}
\textbf{Proof:}
Let $t > 0$ be a parameter to be chosen later. Then
\begin{eqnarray}\label{eq2:thm:Ky Fan norm prob bound for fun of tensors sum}
\mathrm{Pr} \left( \left\Vert g\left( \sum\limits_{j=1}^{m} \mathcal{X}_j  \right)\right\Vert_{(k)} \geq \theta \right)=  \mathrm{Pr} \left(   \left\Vert \exp\left(   t g\left( \sum\limits_{j=1}^{m} \mathcal{X}_j  \right)\right)\right\Vert_{(k)} \geq \exp\left(\theta t \right) \right) \nonumber \\
\leq_1 \exp \left(- \theta t \right) \mathbb{E} \left(   \left\Vert \exp\left(   t g\left( \sum\limits_{j=1}^{m} \mathcal{X}_j  \right)\right)\right\Vert_{(k)} \right) \nonumber \\
\leq_2  \exp \left(- \theta t \right) \mathbb{E} \left(   \left\Vert g \left( \exp\left(   t  \sum\limits_{j=1}^{m} \mathcal{X}_j  \right)\right)\right\Vert_{(k)} \right)
\end{eqnarray}
where $\leq_1$ uses Markov's inequality, $\leq_2$ requires condition provided by Eq.~\eqref{eq:special cond}.

We can further bound the expectation term in Eq.~\eqref{eq1:thm:Ky Fan norm prob bound for fun of tensors sum} as 
\begin{eqnarray}\label{eq3:thm:Ky Fan norm prob bound for fun of tensors sum}
\mathbb{E} \left(   \left\Vert g \left( \exp\left(   t  \sum\limits_{j=1}^{m} \mathcal{X}_j  \right)\right)\right\Vert_{(k)} \right) \leq_3
\mathbb{E} \int_{- \infty}^{\infty} \left\Vert  g\left(   \left\vert  \prod\limits_{j=1}^m e^{ (1 + \iota \tau)t \mathcal{X}_j }   \right\vert \right)  \right\Vert_{(k)} \beta_0(\tau) d \tau \nonumber \\
\leq_4  (n+1)^{s-1}\left( k a_0^s + \sum\limits_{l=1}^{n}a^{ls}_l  \mathbb{E} 
\int_{- \infty}^{\infty} \left\Vert \left\vert \prod\limits_{j=1}^m e^{ (1 + \iota \tau) t \mathcal{X}_j }    \right\vert^{ls} \right\Vert_{(k)}  \beta_0(\tau)  d \tau     \right),~~~
\end{eqnarray}
where $\leq_3$ from Eq.~\eqref{eq2:thm:Multivaraite Tensor Norm Inequalities} in Theorem~\ref{thm:Multivaraite Tensor Norm Inequalities}, $\leq_4$ is obtained from function $g$ definition and Lemma~\ref{lma:Ky Fan Inequalities for the sum of tensosrs}. Again, the expectation term in Eq.~\eqref{eq3:thm:Ky Fan norm prob bound for fun of tensors sum} can be further bounded by Lemma~\ref{lma:Ky Fan Inequalities for the prod of tensosrs} as 
\begin{eqnarray}\label{eq4:thm:Ky Fan norm prob bound for fun of tensors sum}
\mathbb{E} \int_{- \infty}^{\infty} \left\Vert \left\vert  \prod\limits_{j=1}^m e^{ (1 + \iota \tau)  t\mathcal{X}_j }    \right\vert^{ls} \right\Vert_{(k)}  \beta_0(\tau) d \tau & \leq & \mathbb{E} \int_{- \infty}^{\infty}  \sum\limits_{j=1}^m \frac{ \left\Vert \left\vert e^{t \mathcal{X}_j } \right\vert^{p_j ls} \right\Vert_{(k)}  }{p_j}     \beta_0( \tau) d \tau  \nonumber \\
& = & \sum\limits_{j=1}^m \frac{ \mathbb{E} \left\Vert e^{p_j  l s t \mathcal{X}_j } \right\Vert_{(k)} }{p_j}. 
\end{eqnarray}
Note that the final equality is obtained due to that the integrand is indepedent of the variable $\tau$ and $\int_{- \infty}^{\infty} \beta_0(\tau) d \tau= 1$. 

Finally, this theorem is established from Eqs.~\eqref{eq2:thm:Ky Fan norm prob bound for fun of tensors sum},~\eqref{eq3:thm:Ky Fan norm prob bound for fun of tensors sum}, and~\eqref{eq4:thm:Ky Fan norm prob bound for fun of tensors sum}. 
$\hfill \Box$

\textbf{Remarks:} The condition provided by Eq.~\eqref{eq:special cond} can be achieved by normalizing tensors $\mathcal{X}_j$ through scaling.

\subsection{Generalized Tensor Chernoff Bound}\label{sec:Generalized Tensor Chernoff Bound}

Folliowing lemma is about Ky Fan $k$-norm bound for the exponential of a random tensor with the constraint of the maximum singular value. 

\begin{lemma}\label{lma:Ky Fan k norm bound for random tensors max eig}
Given a nonnegative Hermitian random tensor $\mathcal{X} \in \mathbb{C}^{I_1 \times \cdots \times I_N \times I_1 \times \cdots \times I_N}$ with $\lambda_{\max}(\mathcal{X}) \leq \mathrm{R}$.  $ \overline{\frac{ \mathcal{X} +\mathcal{X}^*}{2}}$ is defined as $\mathbb{E}\left( \frac{ \mathcal{X} +\mathcal{X}^*}{2}\right)$, and $ \overline{\frac{ \mathcal{X} - \mathcal{X}^*}{2}}$ is defined as $\mathbb{E}\left( \frac{ \mathcal{X} -\mathcal{X}^*}{2}\right)$, respectively. Then, we have following bound about the expectation value of Ky Fan $k$-norm for the random tensor $\exp( \theta \mathcal{X} )$
\begin{eqnarray}\label{eq2:lma:Ky Fan k norm bound for random tensors  max eig}
\mathbb{E} \left\Vert \exp( \theta \mathcal{X} ) \right\Vert_{(k)} &\leq & k   \left\{ 1  +
 (e^{\theta \mathrm{R}} - 1)    \left[  \sigma_1 \left( \overline{\frac{ \mathcal{X} + \mathcal{X}^*}{2}} \right) +  \sigma_1 \left( \overline{\frac{ \mathcal{X} - \mathcal{X}^*}{2}} \right) \right]+
(e^{\theta \mathrm{R}} - 1 ) C \cdot   \right. ~ \nonumber \\
&   &\left. 
 \left[ \max\limits_i \left( \sum\limits_j \mathbb{E} x^2_{i,j}\right)^{1/2} + 
\max\limits_j \left( \sum\limits_i \mathbb{E} x^2_{i,j}\right)^{1/2} + \left( \sum\limits_{i,j} \mathbb{E} x^4_{i,j}\right)^{1/4} \right.   +\right.  \nonumber \\
&   &\left. \left.
 \max\limits_i \left( \sum\limits_j \mathbb{E} y^2_{i,j}\right)^{1/2} +  \max\limits_j \left( \sum\limits_i \mathbb{E} y^2_{i,j}\right)^{1/2} + \left( \sum\limits_{i,j} \mathbb{E} y^4_{i,j}\right)^{1/4} \right]  \right\}\nonumber \\
&\define&   k   \left\{ 1  +
 (e^{\theta \mathrm{R}} - 1)    \left[  \sigma_1 \left( \overline{\frac{ \mathcal{X} + \mathcal{X}^*}{2}} \right) +  \sigma_1 \left( \overline{\frac{ \mathcal{X} - \mathcal{X}^*}{2}} \right) \right]+ \right. \nonumber \\
&   & \left.  (e^{\theta \mathrm{R}} - 1 ) C \Xi(\mathcal{X}) \right\},
\end{eqnarray}
where $\theta$ is a real number, $C$ is a constant, $x_{i,j}$ and $y_{i,j}$ are entries of matrices  obtained from unfolded random real tensors $\frac{\mathcal{X} + \mathcal{X}^*}{2} -  \overline{\frac{ \mathcal{X} +\mathcal{X}^*}{2}}$ and $\frac{\mathcal{X} - \mathcal{X}^*}{2} -  \overline{\frac{ \mathcal{X} - \mathcal{X}^*}{2}}$, respectively. The matrices from unfolded tensors are obtained by the method presented in Section 2.2~\cite{liang2019further}. For notation simplicity, the term $\Xi(\mathcal{X})$ is defined as 
\begin{eqnarray}\label{eq:abbr of sigma 1 of a random tensor}
\Xi(\mathcal{X})&\define& \left[ \max\limits_i \left( \sum\limits_j \mathbb{E} x^2_{i,j}\right)^{1/2} + 
\max\limits_j \left( \sum\limits_i \mathbb{E} x^2_{i,j}\right)^{1/2} + \left( \sum\limits_{i,j} \mathbb{E} x^4_{i,j}\right)^{1/4} \right.   +  \nonumber \\
&  & \left.
 \max\limits_i \left( \sum\limits_j \mathbb{E} y^2_{i,j}\right)^{1/2} +  \max\limits_j \left( \sum\limits_i \mathbb{E} y^2_{i,j}\right)^{1/2} + \left( \sum\limits_{i,j} \mathbb{E} y^4_{i,j}\right)^{1/4} \right]. 
\end{eqnarray}
\end{lemma}
\textbf{Proof:}
Consider the function $f(x) = e^{\theta x}$. Since $f$ is convex, we have 
\begin{eqnarray}
f(x) \leq f(0) + [f(\mathrm{R}) - f(0)] x~~\mbox{for $x \in [0, \mathrm{R}]$}.
\end{eqnarray}
Because all eigenvalues of $\mathcal{X}$ lie in $[0, \mathrm{R}]$, we will get 
\begin{eqnarray}
e^{\theta \mathcal{X}} \preceq \mathcal{I} + (e^{\theta \mathrm{R}} - 1) \mathcal{X}.
\end{eqnarray}

From majorization inequality for eigenvalues of two positive definite tensors, we have 
\begin{eqnarray}
\mathbb{E} \left\Vert \exp( \theta \mathcal{X} ) \right\Vert_{(k)} &=& 
\sum\limits_{l=1}^{k} \mathbb{E} \sigma_l \left( \exp(\theta \mathcal{X})\right) \nonumber \\
&\leq &  \sum\limits_{l=1}^{k} \mathbb{E} \sigma_l \left( \mathcal{I} +  (e^{\theta \mathrm{R}} - 1) \mathcal{X} \right) \leq k \mathbb{E} \sigma_1 \left( \mathcal{I} +
(e^{\theta \mathrm{R}} - 1) \mathcal{X} \right) 
\end{eqnarray}
where $\sigma_l (\cdot )$ is the $l$-th largest singular value. 

From Theorem G.1.d in~\cite{MR2759813} and Theorem 3.2 in~\cite{liang2019further}, we have $\sigma_1(\mathcal{A} + \mathcal{B}) \leq \sigma_1(\mathcal{A}) + \sigma_1(\mathcal{B})$ for two complex tensors $\mathcal{A}$ and $\mathcal{B}$. Then, we can bound $ \mathbb{E} \sigma_1 \left( \mathcal{I} +
(e^{\theta \mathrm{R}} - 1) \mathcal{X} \right) $ as 
\begin{eqnarray}
 \mathbb{E} \sigma_1 \left( \mathcal{I} +
(e^{\theta \mathrm{R}} - 1) \mathcal{X} \right) & \leq & \sigma_1 \left( \mathcal{I}  \right) + (e^{\theta \mathrm{R}} - 1) \mathbb{E}   \sigma_1(\mathcal{X}) \nonumber \\
&= &1+ (e^{\theta \mathrm{R}} - 1)\mathbb{E}   \sigma_1\left( \frac{ \mathcal{X} +\mathcal{X}^*}{2} + \iota \frac{ \mathcal{X} - \mathcal{X}^*}{2}  \right) \nonumber \\
& \leq &1 + (e^{\theta \mathrm{R}} - 1) \mathbb{E}   \sigma_1\left( \frac{ \mathcal{X} +\mathcal{X}^*}{2} \right) 
+ (e^{\theta \mathrm{R}} - 1) \mathbb{E}  \sigma_1 \left( \iota \frac{ \mathcal{X} - \mathcal{X}^*}{2}  \right)
 \nonumber \\
& \leq &1 + (e^{\theta \mathrm{R}} - 1)    \left[  \sigma_1 \left( \overline{\frac{ \mathcal{X} +\mathcal{X}^*}{2}} \right) +  \mathbb{E} \sigma_1 \left( \frac{ \mathcal{X} +\mathcal{X}^*}{2} -  \overline{\frac{ \mathcal{X} +\mathcal{X}^*}{2}} \right)  \right] +  \nonumber \\
& &
 (e^{\theta \mathrm{R}} - 1)    \left[  \sigma_1 \left( \overline{\frac{ \mathcal{X} -\mathcal{X}^*}{2}} \right) +  \mathbb{E} \sigma_1 \left( \frac{ \mathcal{X} -\mathcal{X}^*}{2} -  \overline{\frac{ \mathcal{X} -\mathcal{X}^*}{2}} \right)  \right]. 
\end{eqnarray}
This lemma is proved by using Theorem 2.5 in~\cite{rudelson2010non} to bound $\mathbb{E}\sigma_1 \left( \frac{ \mathcal{X} +\mathcal{X}^*}{2} -  \overline{\frac{ \mathcal{X} +\mathcal{X}^*}{2}} \right) $ and $ \mathbb{E}\sigma_1 \left(\frac{ \mathcal{X} -\mathcal{X}^*}{2} -  \overline{\frac{ \mathcal{X} -\mathcal{X}^*}{2}} \right) $ since $ \mathbb{E}  \left( \frac{ \mathcal{X} +\mathcal{X}^*}{2} -  \overline{\frac{ \mathcal{X} +\mathcal{X}^*}{2}} \right) =\mathbb{E}  \left( \frac{ \mathcal{X} - \mathcal{X}^*}{2} -  \overline{\frac{ \mathcal{X} - \mathcal{X}^*}{2}} \right) = \mathcal{O}$. 
$\hfill \Box$

Following theorem is presented to provide the general tensor Chernoff bound under Ky Fan $k$-norm. 
\begin{theorem}[Generalized Tensor Chernoff Bound]\label{thm:Generalized Tensor Chernoff Bound}
Consider a sequence $\{ \mathcal{X}_j  \in \mathbb{C}^{I_1 \times \cdots \times I_N  \times I_1 \times \cdots \times I_N} \}$ of independent, random, Hermitian tensors. Let $g$ be a polynomial function with degree $n$ and nonnegative coeffecients $a_0, a_1, \cdots, a_n$ raised by power $s \geq 1$, i.e., $g(x) = \left(a_0 + a_1 x  +\cdots + a_n x^n \right)^s$ with $s \geq 1$. Suppose following condition is satisfied:
\begin{eqnarray}\label{eq:special cond Chernoff Bound 1}
g \left( \exp\left(t \sum\limits_{j=1}^{m} \mathcal{X}_j \right)\right)  \succeq \exp\left(t g \left( \sum\limits_{j=1}^{m} \mathcal{X}_j   \right) \right)~~\mbox{almost surely},
\end{eqnarray}
where $t > 0$. Moreover, we require  
\begin{eqnarray}\label{eq:special cond Chernoff Bound 2}
\mathcal{X}_i \succeq \mathcal{O} \mbox{~~and~~} \lambda_{\max}(\mathcal{X}_i) \leq \mathrm{R}
\mbox{~~ almost surely.}
\end{eqnarray}
Then we have following inequality:
\begin{eqnarray}\label{eq1:thm:GeneralizedTensorChernoffBound}
\mathrm{Pr} \left( \left\Vert g\left( \sum\limits_{j=1}^{m} \mathcal{X}_j  \right)\right\Vert_{(k)}  \geq \theta \right)  \leq  (n+1)^{s-1} \inf\limits_{t > 0} e^{- \theta t } \cdot ~~~~~~~~~~~~~~~~~~~~~~~~~~~~ \nonumber \\
 \left\{ ka_0^s + \sum\limits_{l=1}^{n} \sum\limits_{j=1}^m \frac{k a_l^{ls}}{m} \left[ 1 +\left( e^{mlsRt} - 1 \right) \overline{\sigma_1(\mathcal{X}_j)}  +  C \left( e^{mlsRt} - 1 \right) \Xi(\mathcal{X}_j) \right] \right\},
\end{eqnarray}
where $ \overline{\sigma_1(\mathcal{X}_j)} \define \left[  \sigma_1 \left( \overline{\frac{ \mathcal{X}_j + \mathcal{X}_j^*}{2}} \right) +  \sigma_1 \left( \overline{\frac{ \mathcal{X}_j - \mathcal{X}_j^*}{2}} \right) \right] $. Let us define following three terms $A_1(t), A_2(t)$ and $A_3(t)$ as 
\begin{eqnarray}\label{eq:A 1 2 3 define}
A_1 (t) &=&  ka_0^s + \sum\limits_{l=1}^{n} \sum\limits_{j=1}^m \frac{ k a_l^{ls}}{m} \left[ 1 + \left( e^{mlsRt} - 1 \right) \overline{\sigma_1(\mathcal{X}_j)}+  \right. \nonumber \\
& & \left. C \left( e^{mlsRt} - 1 \right) \Xi(\mathcal{X}_j) \right], \nonumber \\
A_2 (t) &=& \sum\limits_{l=1}^{n} \sum\limits_{j=1}^m klsRk a_l^{ls}  \left( \overline{\sigma_1(\mathcal{X}_j)}  + C \Xi(\mathcal{X}_j) \right) e^{mlsRt}, \nonumber \\
A_3 (t) &=&  \sum\limits_{l=1}^{n} \sum\limits_{j=1}^m km(lsR)^2  a_l^{ls}  \left( \overline{\sigma_1(\mathcal{X}_j)}  + C \Xi(\mathcal{X}_j) \right) e^{mlsRt}. 
\end{eqnarray}
If we have $\theta^2 A_1 (t) - 2 \theta A_2 (t)+ A_3 (t) >0 $ for $t > 0$, then the bound in Eq.~\eqref{eq1:thm:GeneralizedTensorChernoffBound} is a convex function with respect to $t$ and the minimizer, denoted as $t_{opt}$, is the solution of the following equation: 
\begin{eqnarray}
\sum\limits_{l=1}^{n} \sum\limits_{j=1}^m \frac{k a_l^{ls}}{m}e^{mlsRt_{opt}}\left[ (mlsR + \theta) ( \overline{\sigma_1(\mathcal{X}_j)}   + C \Xi(\mathcal{X}_j) )   \right] = \nonumber \\
\theta \left[  ka_0^s + \sum\limits_{l=1}^{n} \sum\limits_{j=1}^m \frac{ k a_l^{ls}}{m}\left( 1 -  \overline{\sigma_1(\mathcal{X}_j)} - C \Xi(\mathcal{X}_j)  \right)\right]
\end{eqnarray}
\end{theorem}
\textbf{Proof:}
From Theorem~\ref{thm:Ky Fan norm prob bound for fun of tensors sum} and Lemma~\ref{lma:Ky Fan k norm bound for random tensors max eig}, we will have the bound given by Eq.~\eqref{eq1:thm:GeneralizedTensorChernoffBound} by taking $p_j =m$. 

The convexity condition, $\theta^2 A_1 (t) - 2 \theta A_2 (t)+ A_3 (t) >0 $, of this generalized Chernoff bound is obtained by setting the second derivative of Eq.~\eqref{eq2:thm:GeneralizedTensorChernoffBound} with respect to $t$ greater than zero. 
\begin{eqnarray}\label{eq2:thm:GeneralizedTensorChernoffBound}
e^{- \theta t } \left\{ ka_0^s + \sum\limits_{l=1}^{n} \sum\limits_{j=1}^m \frac{k a_l^{ls}}{m} \left[ 1 + \left( e^{mlsRt} - 1 \right) \overline{\sigma_1(\mathcal{X}_j)} + \right. \right. \nonumber \\
 \left. \left. C \left( e^{ mlsRt} - 1 \right) \Xi(\mathcal{X}_j) \right] \right\}.~~~~~~~~~~~~~~~~~~~~~~~~~~
\end{eqnarray}
Similarly, the optimizer $t_{opt}$ is obtained by setting the first derivative of Eq.~\eqref{eq2:thm:GeneralizedTensorChernoffBound} with respect to $t$ to zero. 
$\hfill \Box$



\subsection{Generalized Bernstein Tensor Bound}\label{sec:Generalized Bernstein Tensor Bound}


In this section, we will present a generalized tensor Bernstein bound, and we will begin with a lemma to bound exponential of a random tensor.
\begin{lemma}\label{lem:Subexponential Bernstein mgf}
Suppose that $\mathcal{X}$ is a random Hermitian tensor that satisfies
\begin{eqnarray}\label{eq1:lem:Subexponential Bernstein mgf}
\mathcal{X}^p  \preceq \frac{p! \mathcal{A}^2}{2} 
\mbox{~~ almost surely for $p=2,3,4,\cdots$,} 
\end{eqnarray}
where $\mathcal{A}$ is a fixed Hermitian tensor. Then, we have
\begin{eqnarray}
 e^{t \mathcal{X}} \preceq \mathcal{I} + t \mathcal{X} + \frac{t^2 \mathcal{A}^2}{2 ( 1 - t) } 
\mbox{~~ almost surely,} 
\end{eqnarray}
where $0 < t < 1$.
\end{lemma}
\textbf{Proof:}
From Tayler series of the tensor exponential expansion, we have 
\begin{eqnarray}
e^{t \mathcal{X}} &=& \mathcal{I} + t \mathcal{X} + \sum\limits_{p=2}^{\infty} \frac{t^p (\mathcal{X}^p)}{p!} \preceq \mathcal{I}+  t \mathcal{X} +  \sum\limits_{p=2}^{\infty} \frac{t^p \mathcal{A}^2}{2} = \mathcal{I}+  t \mathcal{X}  + \frac{t^2 \mathcal{A}^2}{2 (1-t)}.
\end{eqnarray}
Therefore, this Lemma is proved. 
$\hfill \Box$

Following lemma is about Ky Fan $k$-norm bound for the exponential of a random tensor with subexponential constraints. 
\begin{lemma}\label{lma:Ky Fan k norm bound for random tensors subexp}
Given a Hermitian random tensor $\mathcal{X} \in  \mathbb{C}^{I_1 \times \cdots \times I_N \times I_1 \times \cdots \times I_N}$ with $\mathbb{E} \mathcal{X} = \mathcal{O}$ and 
\begin{eqnarray}\label{eq1:lma:Ky Fan k norm bound for random tensors subexp}
\mathcal{X}^p  \preceq \frac{p! \mathcal{A}^2}{2} 
\mbox{~~ almost surely for $p=2,3,4,\cdots$,} 
\end{eqnarray}
where $\mathcal{A}$ is a positive Hermitian tensor. Then, we have following bound about the expectation value of Ky Fan $k$-norm for the random tensor $\exp( \theta \mathcal{X} )$
\begin{eqnarray}\label{eq2:lma:Ky Fan k norm bound for random tensors  subexp}
\mathbb{E} \left\Vert \exp( \theta \mathcal{X} ) \right\Vert_{(k)} &\leq&  k \left\{1 +  \frac{\theta^2}{ 2 (1 - \theta)} \sigma_1 \left( \mathcal{A}^2\right) +  \right.  \nonumber \\
&  & \left. 
\theta C\left[ \max\limits_i \left( \sum\limits_j \mathbb{E} x^2_{i,j}\right)^{1/2} + 
\max\limits_j \left( \sum\limits_i \mathbb{E} x^2_{i,j}\right)^{1/2} + \left( \sum\limits_{i,j} \mathbb{E} x^4_{i,j}\right)^{1/4} \right]   +\right.  \nonumber \\
&  & \left.
\theta C \left[ \max\limits_i \left( \sum\limits_j \mathbb{E} y^2_{i,j}\right)^{1/2} +  \max\limits_j \left( \sum\limits_i \mathbb{E} y^2_{i,j}\right)^{1/2} + \left( \sum\limits_{i,j} \mathbb{E} y^4_{i,j}\right)^{1/4} \right]  \right\} \nonumber \\
&\define& k \left\{ 1 +  \frac{\theta^2}{ 2 (1 - \theta)} \sigma_1 \left( \mathcal{A}^2\right) +  \theta C \Upsilon(\mathcal{X}) \right\} ,
\end{eqnarray}
where $\theta$ is a positive number in the ange $(0, 1)$, $C$ is a constant, $x_{i,j}$ and $y_{i,j}$ are entries of the unfolded random real tensors $\frac{\mathcal{X} + \mathcal{X}^*}{2}$ and $\frac{\mathcal{X} - \mathcal{X}^*}{2}$, respectively. from Eq.~\eqref{eq:abbr of sigma 1 of a random tensor}.  The matrices from unfolded tensors are obtained by the method presented in Section 2.2~\cite{liang2019further}.  For notation simplicity, the term $\Upsilon$ is defined as 
\begin{eqnarray}\label{eq:abbr of sigma 1 of a random tensor zero mean}
\Upsilon(\mathcal{X})&\define& \left[ \max\limits_i \left( \sum\limits_j \mathbb{E} x^2_{i,j}\right)^{1/2} + 
\max\limits_j \left( \sum\limits_i \mathbb{E} x^2_{i,j}\right)^{1/2} + \left( \sum\limits_{i,j} \mathbb{E} x^4_{i,j}\right)^{1/4} \right.   +  \nonumber \\
&  & \left.
 \max\limits_i \left( \sum\limits_j \mathbb{E} y^2_{i,j}\right)^{1/2} +  \max\limits_j \left( \sum\limits_i \mathbb{E} y^2_{i,j}\right)^{1/2} + \left( \sum\limits_{i,j} \mathbb{E} y^4_{i,j}\right)^{1/4} \right], 
\end{eqnarray}
where random tensor $\mathcal{X}$ has zero tensor as its mean. 
\end{lemma}
\textbf{Proof:}
From Lemma~\ref{lem:Subexponential Bernstein mgf} and majorization inequality for eigenvalues of two positive definite tensors, we have 
\begin{eqnarray}
\mathbb{E} \left\Vert \exp( \theta \mathcal{X} ) \right\Vert_{(k)} &=& 
\sum\limits_{l=1}^{k} \mathbb{E} \sigma_l \left( \exp(\theta \mathcal{X})\right) \nonumber \\
&\leq &  \sum\limits_{l=1}^{k} \mathbb{E} \sigma_l \left( \mathcal{I} + \theta \mathcal{X} + 
 \frac{\theta^2 \mathcal{A}^2}{ 2 (1 - \theta)} \right) \leq k \mathbb{E} \sigma_1 \left( \mathcal{I} + \theta \mathcal{X} + 
 \frac{\theta^2 \mathcal{A}^2}{ 2 (1 - \theta)} \right) 
\end{eqnarray}
where $\sigma_l (\cdot )$ is the $l$-th largest singular value. 

From Theorem G.1.d in~\cite{MR2759813} and Theorem 3.2 in~\cite{liang2019further}, we have $\sigma_1(\mathcal{A} + \mathcal{B}) \leq \sigma_1(\mathcal{A}) + \sigma_1(\mathcal{B})$ for two complex tensors $\mathcal{A}$ and $\mathcal{B}$. Then, we can bound $ \mathbb{E} \sigma_1 \left(\mathcal{I} + \theta \mathcal{X} + 
 \frac{\theta^2 \mathcal{A}^2}{ 2 (1 - \theta)} \right) $ as 
\begin{eqnarray}
 \mathbb{E} \sigma_1 \left(\mathcal{I} + \theta \mathcal{X} + 
 \frac{\theta^2 \mathcal{A}^2}{ 2 (1 - \theta)} \right)& \leq &1 +  \frac{\theta^2}{ 2 (1 - \theta)} \sigma_1 \left( \mathcal{A}^2\right) + \theta \mathbb{E}   \sigma_1(\mathcal{X}) \nonumber \\
&= & 1 +  \frac{\theta^2}{ 2 (1 - \theta)} \sigma_1 \left( \mathcal{A}^2\right) + \theta \mathbb{E}   \sigma_1\left( \frac{ \mathcal{X} +\mathcal{X}^*}{2} + \iota  \frac{ \mathcal{X} - \mathcal{X}^*}{2}  \right) \nonumber \\
& \leq &1 +  \frac{\theta^2}{ 2 (1 - \theta)} \sigma_1 \left( \mathcal{A}^2\right) + \theta \mathbb{E}   \sigma_1\left( \frac{ \mathcal{X} +\mathcal{X}^*}{2} \right) 
+ \theta \mathbb{E}  \sigma_1 \left(  \iota  \frac{ \mathcal{X} - \mathcal{X}^*}{2}  \right)
 \nonumber \\
& = & 1 +  \frac{\theta^2}{ 2 (1 - \theta)} \sigma_1 \left( \mathcal{A}^2\right)+ \theta \mathbb{E}   \sigma_1\left( \frac{ \mathcal{X} +\mathcal{X}^*}{2} \right) 
+ \theta \mathbb{E}  \sigma_1 \left( \frac{ \mathcal{X} - \mathcal{X}^*}{2}  \right)
\end{eqnarray}
This lemma is proved by using Theorem 2.5 in~\cite{rudelson2010non} to bound $ \mathbb{E}   \sigma_1\left( \frac{ \mathcal{X} +\mathcal{X}^*}{2} \right)$ and $\mathbb{E}  \sigma_1 \left( \frac{ \mathcal{X} - \mathcal{X}^*}{2}  \right)$ since $ \mathbb{E}  \left( \frac{ \mathcal{X} +\mathcal{X}^*}{2} \right) = \mathbb{E} \left( \frac{ \mathcal{X} - \mathcal{X}^*}{2}  \right) = \mathcal{O}$. 
$\hfill \Box$

\begin{theorem}[Generalized Tensor Bernstein Bound]\label{thm:Generalized Tensor Bernstein Bound}
Consider a sequence $\{ \mathcal{X}_j  \in \mathbb{C}^{I_1 \times \cdots \times I_N  \times I_1 \times \cdots \times I_N} \}$ of independent, random, Hermitian tensors. Let $g$ be a polynomial function with degree $n$ and nonnegative coeffecients $a_0, a_1, \cdots, a_n$ raised by power $s \geq 1$, i.e., $g(x) = \left(a_0 + a_1 x  +\cdots + a_n x^n \right)^s$ with $s \geq 1$. Suppose following condition is satisfied:
\begin{eqnarray}\label{eq:special cond Bernstein Bound}
g \left( \exp\left(t \sum\limits_{j=1}^{m} \mathcal{X}_j \right)\right)  \succeq \exp\left(t g \left( \sum\limits_{j=1}^{m} \mathcal{X}_j   \right) \right)~~\mbox{almost surely},
\end{eqnarray}
where $t > 0$, and we also have 
\begin{eqnarray}
\mathbb{E} \mathcal{X}_j = \mathcal{O} \mbox{~~and~~} \mathcal{X}^p_j \preceq \frac{p! \mathcal{A}_j^2}{2}
\mbox{~~ almost surely for $p=2,3,4,\cdots$.}
\end{eqnarray}
Then we have following inequality:
\begin{eqnarray}\label{eq1:thm:Generalized Tensor Bernstein Bound}
\mathrm{Pr} \left( \left\Vert g\left( \sum\limits_{j=1}^{m} \mathcal{X}_j  \right)\right\Vert_{(k)}  \geq \theta \right) & \leq & (n+1)^{s-1} \inf\limits_{t > 0} e^{- \theta t }k \cdot \nonumber \\
&  & \left\{a_0^s + \sum\limits_{l=1}^{n} \sum\limits_{j=1}^m  a_l^{ls} \left[ \frac{1}{m} + \frac{ m (lst)^2 \sigma_1(\mathcal{A}_j^2) }{2 (1 - m lst)}+ lstC \Upsilon(\mathcal{X}_j) \right] \right\}.
\end{eqnarray}
Let us define following three terms $B_1, B_2$ and $B_3$ as 
\begin{eqnarray}\label{eq:B 1 2 3 define Bernstein}
B_1 (t) &=&  ka_0^s + \sum\limits_{l=1}^{n} \sum\limits_{j=1}^m k  a_l^{ls} \left[ \frac{1}{m} + \frac{ m (lst)^2 \sigma_1(\mathcal{A}_j^2) }{2 (1 - m lst)}+ lstC \Upsilon(\mathcal{X}_j) \right] , \nonumber \\
B_2 (t)&=&  \sum\limits_{l=1}^{n} \sum\limits_{j=1}^m k a_l^{ls} \left[\frac{ \left( 4ml^2s^2t - 3 l^3s^3m^2 t^2\right) \sigma_1(\mathcal{A}_j^2)    }{2  ( 1 - mlst )^2}+lsC \Upsilon(\mathcal{X}_j) \right], \nonumber \\
B_3 (t) &=& \sum\limits_{l=1}^{n} \sum\limits_{j=1}^m \frac{ k mls  a_l^{ls} \left( 2 ls - ml^2s^2 t \right) \sigma_1(\mathcal{A}_j^2)    }{( 1 -mlst )^3}.
\end{eqnarray}
If we have $\theta ^2 B_1 (t)- 2 \theta B_2(t) + B_2 (t) > 0 $ for $0 < t < \frac{1}{mls}$ and $Cls \Upsilon(\mathcal{X}_j) < \frac{\theta}{m}$ for $1 \leq l \leq n$, and $1 \leq j \leq m$, then the bound in Eq.~\eqref{eq1:thm:Generalized Tensor Bernstein Bound} is a convex function with respect to $t$ and the minimizer, denoted as $t_{opt}$, of this bound is the solution of the following equation
\begin{eqnarray}\label{eq:B 1 2 Bernstein eq}
B_2(t_{opt}) = B_1(t_{opt})\theta. 
\end{eqnarray}
\end{theorem}
\textbf{Proof:}
From Theorem~\ref{thm:Ky Fan norm prob bound for fun of tensors sum} and Lemma~\ref{lma:Ky Fan k norm bound for random tensors subexp}, we will have the bound given by Eq.~\eqref{eq1:thm:Generalized Tensor Bernstein Bound} by taking $p_j =m$. 

The convexity condition of this generalized Bernstein bound and its optimizer $t_{opt}$ are obtained by taking first and second derivatives of Eq.~\eqref{eq2:thm:Generalized Tensor Bernstein Bound} with respect to $t$ through tedious algebraic manipulations. 
\begin{eqnarray}\label{eq2:thm:Generalized Tensor Bernstein Bound}
 e^{- \theta t }k \left\{a_0^s + \sum\limits_{l=1}^{n} \sum\limits_{j=1}^m  a_l^{ls} \left[ \frac{1}{m} + \frac{ m (lst)^2 \sigma_1(\mathcal{A}_j^2) }{2 (1-  mlst)}+ lstC \Upsilon(\mathcal{X}_j) \right] \right\}.
\end{eqnarray}
The conidtion of $\theta ^2 B_1 (t)- 2 \theta B_2(t) + B_2 (t) > 0$ for $0 < t < \frac{1}{mls}$ is obtained by having the second derivative of Eq.~\eqref{eq2:thm:Generalized Tensor Bernstein Bound} is positive. The optimizer $t_{opt}$ is obtained by setting the first derivative of Eq.~\eqref{eq2:thm:Generalized Tensor Bernstein Bound} to zero, which is Eq.~\eqref{eq:B 1 2 Bernstein eq}. The convexity condition of Eq.~\eqref{eq2:thm:Generalized Tensor Bernstein Bound} and $Cls \Upsilon(\mathcal{X}_j) < \frac{\theta}{m}$ for $1 \leq l \leq n$, and $1 \leq j \leq m$ make sure that there is a unique real solution for $t_{opt}$ between $0$ and $\frac{1}{mls}$. 
$\hfill \Box$

\section{Covariance Tensor Characterization by Generalized Tensor Chernoff Inequality}\label{sec:
Covariance Tensor Characterization by Generalized Tensor Chernoff Inequality}

In this section, we will try to apply generalized tensor Chernoff inequality
derived in Section 4.2 to bound Ky Fan norm of covariance tensor. In~\cite{marques2017stationary}, Marques
et al. provide a comprehensive introduction to the spectral analysis and estimation
of graph stationary processes based on graph signal processing (GSP). We extend
their settings from vectors/matrices used in traditional GSP to hypergraph signal
processing, where tensors are applied to characterize high-dimensional signals~\cite{zhang2019introducing}.

Let $\mathfrak{G}= (\mathfrak{N}, \mathfrak{E})$ be a directed hypergraph with nodes set $\mathfrak{N}$ and directed edges set $\mathfrak{E}$ such that if there exists a hyperedge between two sets of $M$ nodes $(i_1, \cdots, i_M,$ $j_1, \cdots, j_M) \in \mathfrak{E}$. We associate $\mathfrak{G}$ with the hypergraph shift operator (HGSO) $\mathcal{S}$, defined as an square tensor with dimensions $I_1 \times \cdots \times I_M \times I_1 \times \cdots \times I_M$ whose entry $s_{(i_1, \cdots, i_M, j_1, \cdots, j_M)} \neq 0$ if $(i_1, \cdots, i_M, j_1, \cdots, j_M) \in \mathfrak{E}$. We introduce a hypergraph filter $\mathcal{H} : \mathbb{C}^{I_1 \times \cdots \times I_M} \rightarrow  \mathbb{C}^{I_1 \times \cdots \times I_M}$, defined as a linear hypergraph signal operator with the form 
\begin{eqnarray}\label{eq:operator H}
\mathcal{H} \define \sum\limits_{k=0}^{K-1}h_k \mathcal{S}^k,
\end{eqnarray}
where $h_k$ are scaler coefficients. The covariance tensor of output signals $x$ after
ltering white input signals by hypergraph filter shown in Eq.~\eqref{eq:operator H} will be expressed as
\begin{eqnarray}\label{eq:cov ten}
\mathcal{C}_x(\mathbf{h}) &=& \mathcal{H}^{H} \star_M \mathcal{H} = \sum\limits_{k=0, k'=0}^{K-1} h_k h_{k'} \mathcal{S}^k \star_M \left( \mathcal{S}^H \right)^k  \nonumber \\
&=_1& \sum\limits_{k=0}^{2(K-1)} \gamma_k \mathcal{S}^k,
\end{eqnarray}
where $=_1$ is true if HGSO $\mathcal{S}$ is a symmetric tensor, i.e., $s_{i_1, \cdots, i_M, j_1, \cdots, j_M} =$  \\ $s_{ j_1, \cdots, j_M, i_1, \cdots, i_M}$. The coefficients $\gamma_k \define \sum\limits_{k' + k'' = k} h_{k'} h_{k''}$. 

It is shown by the work~\cite{navarro2020joint} that although the correlation information of signal is
given by the \emph{dense} tensor, the actual relation is easier to be described by the more
sparse tensor $\mathcal{S}$. Examples about relationships between the HGSO and the covariance
tensor $\mathcal{C}_x(\mathbf{h})$ include
\begin{itemize}
\item $\mathcal{C}_x(\mathbf{h}) = \sum\limits_{k=0}^{2(K-1)} \gamma_k \mathcal{S}^k$, as in graph filtering;
\item $\mathcal{C}_x(\mathbf{h}) = \mathcal{S}^{-1}$, as in in conditionally independent Markov random fields;
\item  $\mathcal{C}_x(\mathbf{h}) = (\mathcal{I} - \mathcal{S})^{-2}$, as in symmetric structural equation models with white exogenous inputs.
\end{itemize}

In the sequel, we will bound the Ky Fan norm for the covariance tensor $\mathcal{C}_x(\mathbf{h})$ when $\mathbf{h} = [h_0, h_1]$. In random environment, suppose HGSO $\mathcal{S}$ is obtained by sample average as
\begin{eqnarray}
\mathcal{S}&=& \frac{1}{m} \sum\limits_{j=1}^m \mathcal{X}_j =  \sum\limits_{j=1}^m \mathcal{X}'_j ,
\end{eqnarray}
where $\mathcal{X}' =\frac{\mathcal{X}}{m}$. Since the graph filter coefficients are  $\mathbf{h} = [h_0, h_1]$, from Eq.~\eqref{eq:cov ten}, the corresponding polynomial relation between $\mathcal{C}_x(\mathbf{h}) $ and $\mathcal{S}$ is
\begin{eqnarray}
\mathcal{C}_x([h_0, h_1]) = h_0^2 + 2 h_0 h_1 \mathcal{S} + h_1^2 \mathcal{S}^2,
\end{eqnarray}
which is the polynomial function $g(x) = (a_0+a_1 x + a_2x^2)^1 = h_0^2 + 2 h_0 h_1 x + h^2_1 x^2$ in Theorem~\ref{thm:Generalized Tensor Chernoff Bound}. We assume that random sampled tensors $\mathcal{X}'_j$ are identical distributed as $\mathcal{X}'$ are satisfy Eq.~\eqref{eq:special cond Chernoff Bound 1} and Eq.~\eqref{eq:special cond Chernoff Bound 2}. Then we have following bound of Ky Fan norm for the covaraince $\mathcal{C}_x[h_0, h_1]$ from Theorem~\ref{thm:Generalized Tensor Chernoff Bound}:
\begin{eqnarray}
\mathrm{Pr}\left( \left\Vert \mathcal{C}_x ( [h_0, h_1] )  \right\Vert_{(k)} \geq \theta \right)
\leq  \inf\limits_{t > 0} ke^{- \theta t} ~~~~~~~~~~~~~~~~~~~~~~~~~~~~~~~~~~~~~~~~\nonumber \\
       \cdot \left\{ h_0^2 + \sum\limits_{l=1}^{2}a_l^l \left[ 1 + ( e^{mlRt}- 1) \overline{\sigma_1 (\mathcal{X}')} + C(e^{mlRt} - 1) \Xi(\mathcal{X}')    \right] \right\},
\end{eqnarray}
where $\overline{\sigma_1 (\mathcal{X}')} \define \left[\sigma_1\left(\overline{\frac{\mathcal{X}' + \mathcal{X}^{'*}}{2}} \right) + \sigma_1\left( \overline{\frac{\mathcal{X}' - \mathcal{X}^{'*}}{2}}  \right) \right]$ and $\Xi(\mathcal{X})$ is defined as~\eqref{eq:abbr of sigma 1 of a random tensor}.

\section{Conclusions}\label{sec:Conclusions}

In this work, we generalize previous work by considering the tail bound behavior of the top $k$-largest singular values of a function of random tensors summation. In previous work, we only considered the tail bound behavior of the largest singular value of tensors summation directly (identity function). Majorization and antisymmetric tensor products are our main gadgets used to derive bounds for unitarily norms of multivariate tensors. Then, we apply Laplace transform method to these bounds of unitarily norms of multivariate tensors to obtain Ky Fan $k$-norm concentration inequalities for a function of tensors summation. Under this approach, generalized tensor Chernoff and Bernstein inequalities are special cases of Ky Fan $k$-norm concentration inequalities obtained by restricting different random tensors conditions.

Possible future works include considering a more general unitarily invariant norm instead of Ky Fan $k$-norm in our tail bounds and other functions of random tensors summation besides the power of polynomials. 

\section*{Acknowledgments}
The helpful comments of the referees are gratefully acknowledged.

\bibliographystyle{siamplain}
\bibliography{GenTensorConcenIndep_Bib}
\end{document}